\documentclass[11pt]{article}

\usepackage{amssymb}
\usepackage{amsmath}
\usepackage{amscd}
\usepackage[dvips]{graphicx}
\usepackage{mathrsfs}
\usepackage{wrapfig} 
\usepackage{picins}
\usepackage{picinpar}
\usepackage{psfig}
\usepackage{rotating}

\usepackage{setspace}
\usepackage{float}
\restylefloat{table}

\makeatletter
\renewcommand\@biblabel[1]{#1.}
\makeatother

\newcommand{\tabincell}[2]{\begin{tabular}{@{}#1@{}}#2\end{tabular}}

\begin{document}

\title{Asymptotic Analysis of Spectral Properties of Finite Capacity Processor Shared Queues}

\author{
{\sc Qiang Zhen}\thanks{
Address for correspondence: Department of Mathematics and Statistics,
University of North Florida, 1 UNF Dr, Bldg. 14/2716, Jacksonville, FL 32224, USA.
{\em Email:} q.zhen@unf.edu.}
\and
and
\and
{\sc Charles Knessl}\thanks{
C. Knessl was partly supported by NSA grants H 98230-08-1-0102 and H 98230-11-1-0184. 
}}
\date{January 23, 2013 }
\maketitle

  \begin{center}
    \rule[1ex]{1\textwidth}{.5pt}
  \end{center}
\noindent We consider sojourn (or response) times in processor-shared queues that have a finite customer capacity. Computing the response time of a tagged customer involves solving a finite system of linear ODEs. Writing the system in matrix form, we study the eigenvectors and eigenvalues in the limit as the size of the matrix becomes large. This corresponds to finite capacity models where the system can only hold a large number $K$ of customers. Using asymptotic methods we reduce the eigenvalue problem to that of a standard differential equation, such as the Airy equation. The dominant eigenvalue leads to the tail of a customer's sojourn time distribution. Some numerical results are given to assess the accuracy of the asymptotic results.  
  \begin{center}
    \rule[1ex]{1\textwidth}{.5pt}
  \end{center}

\section{Introduction}\label{section:1}
The study of processor shared queues has received much attention over the past 45 or so years. The processor sharing (PS) discipline has the advantage over, say, first-in first-out (FIFO), in that shorter jobs tend to get through the system more rapidly. In recent years there has been renewed attention paid to such models, due to their applicability to the flow-level performance of bandwidth-sharing protocols in packet-switched communication networks (see \cite{HEY}-\cite{NAB}). 

Perhaps the simplest example of such a model is the $M/M/1$-PS queue. Here customers arrive according to a Poisson process with rate parameter $\lambda$, the server works at rate $\mu$, there is no queue, and if there are $\mathcal{N}(t)(>0)$ customers in the system each gets an equal fraction $(=1/\mathcal{N}(t))$ of the server. Setting the traffic intensity $\rho$ as $\rho=\lambda/\mu$ it is well known that the steady state distribution of $\mathcal{N}(t)$ is the geometric distribution $\Pr[\mathcal{N}(\infty)=n]=(1-\rho)\rho^n$, which exists for $\rho<1$, and this result is the same as that for the standard FIFO $M/M/1$ model. However, PS and FIFO models differ significantly if we consider the ``sojourn time''. This is defined as the time it takes for a given customer, called a ``tagged customer'', to get through the system (after having obtained the required amount of service). The sojourn time is a random variable that we denote by $\mathcal{V}$. For the simplest $M/M/1$ model, the distribution of $\mathcal{V}$ depends on the total service time $\mathcal{X}$ that the customer requests and also on the number of other customers present when the tagged customer enters the system.

There are two natural variants of the $M/M/1$-PS model that put an upper bound on the number of customers that can be served by the processor. These are the finite population model and the finite capacity model. In the finite population model there are a total of $N$ customers, and each customer will enter service in the next $\Delta t$ time units with probability $\lambda_0\Delta t+o(\Delta t)$. At any time there are $\mathcal{N}(t)$ customers being served and the remaining $N-\mathcal{N}(t)$ customers are in the general population. Hence the total arrival rate is $\lambda_0[N-\mathcal{N}(t)]$ and we may view the model as a PS queue with a state-dependent arrival rate that decreases linearly to zero. Once a customer finishes service that customer re-enters the general population. The service times are exponentially distributed with mean $1/\mu$ and we define the traffic intensity $\rho$ by $\rho=\lambda_0N/\mu$. This model may describe, for example, a network of $N$ terminals in series with a processor-shared CPU. This may be viewed as a closed two node queueing network. 

The finite population model does not see amenable to an exact solution. However, various asymptotic studies have been done in the limit $N\to\infty$, so that the total population, or the number of terminals, is large (see \cite{MIT}--\cite{MOR_C88}). In \cite{ZHE_O12} we studied the spectral structure of the finite population model as $N\to\infty$, with three cases of $\rho$: $\rho<1$, $\rho-1=O(N^{-2/3})$ and $\rho>1$.

A second variant is the finite capacity model where the processor can serve at most $K$ customers. Thus if $\mathcal{N}(t)=K$ and a further arrival occurs, that customer is turned away and lost. This can also be viewed as a PS queue with a state-dependent arrival rate $\lambda(\mathcal{N}(t))$, with $\lambda(\mathcal{N}(t))=\lambda$ if $\mathcal{N}(t)<K$ and $\lambda(\mathcal{N}(t))=0$ if $\mathcal{N}(t)=K$. For this model we set $\rho=\lambda/\mu$, as in the infinite capacity case where $K=\infty$.

Some previous work on the finite capacity model appears in \cite{KNE_O90}-\cite{ZHE_O09}. In \cite{KNE_O90} we analyzed the conditional moments $E[\mathcal{V}^l|\mathcal{X}=x,\mathcal{N}(0^-)=n]$ and in \cite{KNE_O93} the unconditional sojourn time distribution $p(t)$, in the limit of large capacities $K$. The asymptotics tend to be very different for $\rho<1$, $\rho>1$, $\rho-1=O(K^{-1})$ and $\rho-1=O(K^{-1/2})$. More recently, in \cite{ZHE_O09} we gave an explicit, albeit complicated, exact expression for the conditional sojourn time density in the finite capacity model. However, evaluating this exact solution requires an inverse Laplace transform that corresponds to a contour integral with an integrand containing special functions related to hypergeometric functions. Due to this complexity we instead use a singular perturbation approach to study the spectral properties of the finite capacity model. Such an approach should also be applicable to models with other state-dependent arrival rates, such as the finite population model in \cite{ZHE_O12} and queues with discouraged arrival and balking. The analysis here is completed, independent of that in \cite{ZHE_O09}.

We also mention some work on a related finite capacity PS model, where the processor can also serve at most $K$ customers, but now if the system is filled to capacity further arrivals are not lost but rather placed in waiting positions. When the number of customers being served reaches $K-1$ via a departure, one of the waiting customers is placed into service. Various rules as to which customer is allowed to enter service lead to different variants of this model. Some exact, but again very complicated, expressions for these models appear in \cite{AVI_R} and \cite{AVI_S}, while heavy traffic approximations are obtained in \cite{ZHA}.

In this paper we study the spectral structure of the finite capacity model as $K\to\infty$. We denote the sojourn time by $\mathcal{V}=\mathcal{V}(K)$ and its conditional density we call $p_n(t)$ with 
\begin{equation}\label{s1_pnt}
p_n(t)dt=\Pr\Big[\mathcal{V}(K)\in(t,t+dt)\Big|\mathcal{N}(0^-)=n\Big].
\end{equation}
Here $\mathcal{N}(0^-)$ denotes the number of other customers present in the system immediately before the tagged customer arrives, and thus $0\le \mathcal{N}(0^-)\le K-1$. Then we define the column vector $\mathbf{p}(t)=(p_{_0}(t),p_{_1}(t),...,p_{_{K-1}}(t))^T$. $\mathbf{p}(t)$ satisfies a system of ODEs in the form $\mathbf{p}'(t)=\mathbf{B}\mathbf{p}(t)$ where $\mathbf{B}$ is a $K\times K$ tridiagonal matrix, whose entries depend on $\rho=\lambda/\mu$ and $K$. Then eigenvalues of $\mathbf{B}$ are all negative and we denote them by $-\nu_j$ $(j=0,1,..., K-1)$ with the corresponding eigenvectors being $\phi_j(n)=\phi_j(n;K,\rho)$. We shall study this eigenvalue problem for $K\to\infty$ and three cases of $\rho$: $\rho<1$, $\rho>1$ and $\rho-1=O(K^{-2/3})$. In each case we obtain expansions of the $\nu_j$ and then the $\phi_j(n)$, for various ranges of $n$. Often the eigenvectors can be expressed in terms of Airy functions for $K\to\infty$. Since $\mathbf{B}$ is a finite matrix the spectrum is purely discrete, but as the size of the matrix becomes large we sometimes see the eigenvalues coalescing about a certain value. Ordering the eigenvalues as $\nu_0<\nu_1<\nu_2<...$, the tail behavior of $p_n(t)$ and $p(t)$ for $t\to\infty$ is determined by the smallest eigenvalue $\nu_0$. Here $p(t)$ is the unconditional sojourn time density, with 
\begin{equation}\label{s1_pt}
p(t)=\sum_{n=0}^{K-1}p_n(t)\mathrm{Pr}\big[\mathcal{N}(0^-)=n\big]=\sum_{n=0}^{K-1}p_n(t)\frac{(1-\rho)\rho^n}{1-\rho^K},
\end{equation}
as $\mathrm{Pr}\big[\mathcal{N}(0^-)=n\big]=(1-\rho)\rho^n/(1-\rho^K)$, which coincides with the steady state probability of finding $n$ customers in a finite capacity $M/M/1$ queue with capacity $K-1$.

We shall show that the analysis for the finite capacity model is much different than that of the finite population model \cite{ZHE_O12}, and the spectrum will now, for $K\to\infty$, involve Airy functions rather than Hermite polynomials. We shall also see that the zeros of $\phi_j(n)$ now will tend to be concentrated in ranges where $n/K\approx 1$, whereas for the finite population model this concentration occurred where $n/N\approx 0\; (\rho<1,\rho\sim 1)$ or $n/N\approx 1-\rho^{-1}\; (\rho>1)$. Since $\phi_j(n)$ are functions of the discrete variable $n$, by ``zeros" we refer to sign changes of the eigenvectors.

Our basic approach is to use singular perturbation methods to analyze the system of ODEs when $K$ becomes large. The problem can then be reduced to solving simpler, single differential equations whose solutions are known, such as Airy equations. Our analysis does make some assumptions about the forms of various asymptotic series, and about the asymptotic matching of expansions on different scales. We also comment that we assume that the eigenvalue index $j$ is $O(1)$; thus we are not computing the large eigenvalues here. Obtaining the large eigenvalues would likely need a very different asymptotic analysis.

This paper is organized as follows. In section \ref{section:2} we state the mathematical problem and obtain the basic equations. In section \ref{section:3} we summarize our final asymptotic results for the eigenvalues and the (unnormalized) eigenvectors, as well as the tail behaviors of the unconditional density. The derivations are relegated to section \ref{section:4}. Some numerical studies appear in section \ref{section:5}; these assess the accuracy of the asymtptotics.

\section{Statement of the problem}\label{section:2}

We consider the finite capacity $M/M/1/K$-PS model with arrival rate $\lambda$, service rate $\mu=1$, and capacity $K$. Thus the traffic intensity is $\rho=\lambda/\mu=\lambda$. Then the conditional sojourn time density $p_n(t)$ in (\ref{s1_pnt}) satisfies the differential equations 
\begin{equation}\label{s2_fs_pn't}
p_n'(t)=\rho\, p_{n+1}(t)+\frac{n}{n+1}p_{n-1}(t)-(1+\rho)p_n(t),\; 0\le n\le K-2
\end{equation}
with the boundary equations
\begin{equation}\label{s2_fs_pn'tb1}
p_0'(t)=\rho\, p_1(t)-(1+\rho)p_0(t),
\end{equation}
\begin{equation}\label{s2_fs_pn'tb2}
p_{_{K-1}}'(t)=\frac{K-1}{K}p_{_{K-2}}(t)-p_{_{K-1}}(t),
\end{equation}
and the initial condition $p_n(0)=1/(n+1)$. Equations (\ref{s2_fs_pn't})--(\ref{s2_fs_pn'tb2}) are derived in more detail in \cite{KNE_O93} and \cite{COF}. Note that the coefficient $n/(n+1)$ that multiplies $p_{n-1}(t)$ in (\ref{s2_fs_pn't}) corresponds to a departure of a customer other than the tagged one, whose sojourn time we are calculating. If we introduce $p_{_K}(t)$ by requiring (\ref{s2_fs_pn't}) to hold also at $n=K-1$, then (\ref{s2_fs_pn't}) with $n=K-1$ and (\ref{s2_fs_pn'tb2}) leads to 
\begin{equation}\label{s2_fs_pkt}
p_{_K}(t)=p_{_{K-1}}(t).
\end{equation}
This gives an ``artificial boundary condition'' whose use simplifies some of the calculations; the function $p_{_K}(t)$ has no probabilistic meaning. 

Now (\ref{s2_fs_pn't})--(\ref{s2_fs_pn'tb2}) is a finite system of linear constant-coefficient ODEs, which can be written in terms of the tridiagonal matrix $\mathbf{B}$ given in (\ref{s2_matrixB}). Hence, the solution to (\ref{s2_fs_pn't})-(\ref{s2_fs_pn'tb2}) is thus given by the spectral expansion 
\begin{equation}\label{s2_fs_pntspe}
p_n(t)=\sum_{j=0}^{K-1}e^{-\nu_j(K,\rho)t}c_j\phi_j(n;K,\rho)
\end{equation}
where the coefficients $c_j$ follow from 
\begin{equation}\label{s2_fs_phi}
\frac{1}{n+1}=\sum_{j=0}^{K-1}c_j\phi_j(n),\; 0\le n\le K-1.
\end{equation}
Using the orthogonality of the $\phi_j(n)$ leads to the explicit expression 
\begin{equation}\label{s2_cj}
c_j=\frac{\sum_{n=0}^{K-1}\rho^n\,\phi_j(n)}{\sum_{n=0}^{K-1}\rho^n(n+1)\,\phi_j^2(n)}.
\end{equation}
In view of (\ref{s2_fs_pn't})--(\ref{s2_fs_pn'tb2}), $-\nu_j$ are the eigenvalues of the $K\times K$ matrix $\mathbf{B}=\mathbf{B}(K,\rho)$ with
\begin{equation}\label{s2_matrixB}
\mathbf{B}=\left[\begin{array}{ccccccccc}
-1-\rho & \rho & 0 & 0 & 0 & \cdots & 0 & 0 & 0 \\
{1}/{2} & -1-\rho & \rho & 0 & 0 & \cdots & 0 & 0 & 0 \\
0 & 2/3 & -1-\rho & \rho & 0 & \cdots & 0 & 0 & 0 \\
\vdots  & \vdots  & \vdots  & \vdots  & \vdots & \ddots & \vdots  & \vdots & \vdots  \\
0 & 0 & 0 &0 & 0 & \cdots & \frac{K-2}{K-1}  & -1-\rho & \rho  \\
0 & 0 & 0 &0 & 0 & \cdots & 0  & \frac{K-1}{K} & -1
\end{array}\right].
\end{equation}
We shall analyze the $\nu_j$ and $\phi_j(n)$ for $K\to\infty$ with the eigenvalue index $j=O(1)$. We note that for times sufficiently large and $n=O(1)$, the tail of the sojourn time is given by 
\begin{equation*}
p_n(t)\sim c_0\phi_0(n)e^{-\nu_0t},\; t\to\infty.
\end{equation*}
We shall show that the boundary condition in (\ref{s2_fs_pn'tb2}), which corresponds to the last row of the matrix $\mathbf{B}$, will be very important to the structure of the eigenvalue problem. In fact the zeros of the eigenvectors will mostly be concentrated in the range where $n/K\approx 1$. We note that this will be in sharp contrast to the structure of the matrix $\mathbf{A}$ for the finite population model (see \cite{ZHE_O12}, Section 2), where the last row of $\mathbf{A}$ did not effect the eigenvalues asymptotically for $N\to\infty$. 

The unconditional sojourn time density $p(t)$, defined in (\ref{s1_pt}), is asymptotically
\begin{equation}\label{s2_fc_ptsim}
p(t)\sim\bigg[c_0\sum_{n=0}^{K-1}\phi_0(n)\frac{(1-\rho)\rho^n}{1-\rho^K}\bigg]e^{-\nu_0t},\; t\to\infty.
\end{equation}
We will show in section \ref{section:3.5} that for $K\to\infty$ the analysis of $p(t)$ is very different for $\rho<1$, $\rho>1$, and $\rho\sim 1$. 

As a cautionary note, we comment that the analysis of the $\nu_j$ and $\phi_j(n)$ can only be used to calculate $p_n(t)$ and $p(t)$ for very large times $t$. For $K$ large, the structure of the conditional sojourn time density is different not only for different ranges of $\rho$, but also for different space/time scales, where $n$ and $t$ are scaled using $K$. Here we do not attempt to identify precisely how large time must be, relative to $K$, for (\ref{s2_fc_ptsim}) to hold.

\section{Summary of results}\label{section:3}

We summarize our final results for $\nu_j=\nu_j(K,\rho)$ and $\phi_j(n)=\phi_j(n;K,\rho)$. We shall consider separately the cases $\rho<1$, $\rho>1$ and $\rho-1=O(K^{-2/3})$. For each range of $\rho$ we first give the expansion of the eigenvalues as $K\to\infty$, then give the eigenvectors for that spatial range where the sign changes of $\phi_j(n)$ are concentrated, and finally for other spatial ranges of $\xi=n/K$. We comment that the expansion of the first eigenvalue $\nu_0=\nu_0(K,\rho)$ was already obtained in \cite{ZHE_O09}, for all three cases of $\rho$.

\subsection{The case $\rho<1$}\label{section:3.1}

When $\rho<1$ we shall show that
\begin{equation}\label{s32_nu}
\nu_j=(1-\sqrt{\rho})^2+\frac{\sqrt{\rho}}{K}-\frac{\sqrt{\rho}}{K^{4/3}}r_j+\frac{8\sqrt{\rho}}{15K^{5/3}}r_j^2+O(K^{-2}),\; j\ge 0,
\end{equation}
where $r_j$ are the roots of the Airy function $\mathrm{Ai}(\cdot)$, hence $\mathrm{Ai}(r_j)=0$ for $j=0,1,2,...$. We order the roots as $|r_0|<|r_1|<|r_2|<\cdots$, and it is well known that $r_j<0$ and $r_0\approx -2.338$. 

We see that the dependence on $j$ occurs only in the third term in the asymptotic series in (\ref{s32_nu}), and as in the finite population model (\cite{ZHE_O12})  the eigenvalues coalesce about the $M/M/1$ queue relaxation rate $(1-\sqrt{\rho})^2$. However, the expansion is now in powers of $K^{-1/3}$, whereas the finite population model involved powers of $N^{-1/4}$.

To give the eigenvectors when $\rho<1$ we first consider the scale $n=K-O(K^{2/3})$ and introduce $S$ by
\begin{equation*}\label{s32_S}
n=K-K^{2/3}S,\quad \textrm{i.e.,  }S=\frac{K-n}{K^{2/3}}.
\end{equation*} 
Then on the $S$-scale,
\begin{equation}\label{s32_phiS}
\phi_j(n)\sim {k_0}\rho^{-n/2}\mathrm{Ai}(S+r_j),\; j\ge 0
\end{equation}
where ${k_0}=k_0(j)$ is a generic normalization constant. Note that for $j=0$, $\mathrm{Ai}(S+r_0)$ is strictly positive for $S>0$, vanishing only at $S=0$, while $\mathrm{Ai}(S+r_j)$ has $j$ zeros for $S>0$, in addition to vanishing at $S=0$. We also obtain the correction term to (\ref{s32_phiS}), which is given in subsection \ref{section:4.1} (see (\ref{s51_Phiexpand}) and (\ref{s51_Phi(1)})). The expansion in (\ref{s32_phiS}) breaks down for small $S$, since the leading term vanishes.

For the scale $n=K-O(1)$ we introduce $l$ by
\begin{equation*}
n=K-l,\quad \textrm{i.e.,  }l=K-n
\end{equation*} 
and then the eigenvectors behave as
\begin{equation}\label{s32_phil}
\phi_j(n)\sim {k_0}\,\rho^{-n/2}K^{-2/3}\mathrm{Ai}'(r_j)\Big(l+\frac{\sqrt{\rho}}{1-\sqrt{\rho}}\Big),\; j\ge 0.
\end{equation}
In view of (\ref{s32_phiS}) and (\ref{s32_phil}) we conclude that $\phi_0(n)$ has no sign changes (at least for $n=K-O(K^{2/3})$), while $\phi_j(n)$ has exactly $j$ sign changes in $n$, and these are spaced by $O(K^{2/3})$. 

Different expansions must also be constructed for $n=K\xi=O(K)$ with $0<\xi<1$ and for $n=O(1)$. On the $\xi$-scale we obtain
\begin{equation}\label{s32_phixi}
\phi_j(n)\sim{k_1}\rho^{-n/2}K^{-1/12}L(\xi)\exp\Big\{\sqrt{K}\psi(\xi)+K^{1/6}\psi^{(1)}_j(\xi)\Big\}
\end{equation}
where
\begin{equation}\label{s32_psixi}
\psi(\xi)=\sqrt{\xi(1-\xi)}-\sin^{-1}\big(\sqrt{1-\xi}\big),
\end{equation}
\begin{equation}\label{s32_psixi1}
\psi^{(1)}_j(\xi)=-\frac{1}{2}r_j\Big[\sqrt{\xi(1-\xi)}+\sin^{-1}\big(\sqrt{1-\xi}\big)\Big],
\end{equation}
\begin{equation}\label{s32_lxi}
L(\xi)=\Big[\xi(1-\xi)\Big]^{-1/4}.
\end{equation}
By asymptotically matching (\ref{s32_phixi}) for $\xi \uparrow 1$ to (\ref{s32_phiS}) for $S\to\infty$ we can relate the constants ${k_0}$ and ${k_1}$ as follows:
\begin{equation}\label{s32_k0andk1}
{k_1}=\frac{1}{2\sqrt{\pi}}{k_0}.
\end{equation}
Finally, for $n=O(1)$ we obtain
\begin{equation}\label{s32_phio1}
\phi_j(n)\sim {k_2}\rho^{-n/2}\frac{1}{2\pi i}\oint\frac{1}{z^{n+1}}\frac{1}{1-z}\exp\Big(\frac{1}{1-z}\Big)dz,
\end{equation}
where the contour integral is a small loop about $z=0$.
By asymptotically matching (\ref{s32_phixi}) as $\xi\downarrow 0$ to (\ref{s32_phio1}) as $n\to\infty$ we find that 
\begin{equation}\label{s32_k1andk2}
{k_2}={k_1}\frac{2\sqrt{\pi}}{\sqrt{e}}K^{1/6}\exp\Big\{-\frac{\pi}{2}\sqrt{K}-\frac{\pi}{4}r_jK^{1/6}\Big\}.
\end{equation}

\subsection{The case $\rho>1$}\label{section:3.2}

Now consider the case $\rho>1$ with $K\to\infty$. The eigenvalues and eigenfunctions now behave very differently according as $j=0$ or $j\ge 1$. The first eigenvalue is $O(K^{-1})$ as $K\to\infty$ and has the expansion 
\begin{equation}\label{s32_nu0}
\nu_0=\frac{1}{K}+\frac{1}{\rho-1}\frac{1}{K^2}+\frac{1}{(\rho-1)^2}\frac{1}{K^3}+O(K^{-4}).
\end{equation}
The corresponding eigenvector has the following behaviors on the $\xi$ and $n$-scales: 
\begin{equation}\label{s32_phi0xi}
\phi_0(n)\sim k_0^*\,\xi^{\frac{1}{\rho-1}}\exp\Big(-\frac{\xi}{\rho-1}\Big),\; 0<\xi\le 1,
\end{equation}
\begin{equation}\label{s32_phi0n}
\phi_0(n)\sim k_1^*\frac{e^{i\pi\rho/(\rho-1)}}{2\pi i}\int_\mathcal{C}(1-z)^{-\frac{\rho}{\rho-1}}\Big(z-\frac{1}{\rho}\Big)^{\frac{1}{\rho-1}}z^ndz,
\end{equation}
where $\mathcal{C}$ is a closed loop that encircles the branch cut, where $\Im(z)=0$ and $\Re(z)\in[\rho^{-1},1]$, in the $z$-plane, with the integrand being analytic exterior to this cut. We also obtain the correction term to (\ref{s32_phi0xi}), which is given in subsection \ref{section:4.2}. The constants $k_0^*$ and $k_1^*$ are related by
\begin{equation}\label{s32_k0*andk1*}
k_0^*=\frac{K^{\frac{1}{\rho-1}}}{\Gamma\big(\frac{\rho}{\rho-1}\big)}\Big(1-\frac{1}{\rho}\Big)^{\frac{1}{\rho-1}}k_1^*.
\end{equation}
Thus the smallest eigenvalue $\nu_0$ has an eigenfunction which is ``spread out'' over the entire interval $\xi\in(0,1)$, with a distortion for small values of $\xi=n/K$, where (\ref{s32_phi0n}) applies.

The other eigenvalues $\nu_j$, $j\ge 1$ are similar in form to (\ref{s32_nu}), in that 
\begin{equation}\label{s32_nuj}
\nu_j=(\sqrt{\rho}-1)^2+\frac{\sqrt{\rho}}{K}-\frac{\sqrt{\rho}}{K^{4/3}}r_{j-1}+\frac{8\sqrt{\rho}}{15K^{5/3}}r_{j-1}^2+O(K^{-2}),\; j\ge 1
\end{equation}
where $r_{j-1}$ are again the roots of the Airy function. Thus for $j\ge 1$ the eigenvalues are $O(1)$ and coalesce near $(\sqrt{\rho}-1)^2$. The eigenfunctions for $j\ge 1$ are given by
\begin{equation}\label{s32_phiS2}
\phi_j(n)\sim {k_0}\rho^{-n/2}\mathrm{Ai}(S+r_{j-1}),\; j\ge 1,
\end{equation}
\begin{equation}\label{s32_phil2}
\phi_j(n)\sim {k_0}\rho^{-n/2}K^{-2/3}\mathrm{Ai}'(r_{j-1})\Big(l-\frac{\sqrt{\rho}}{\sqrt{\rho}-1}\Big),\; j\ge 1
\end{equation}
on the $S$ and $l$ scales, respectively. Note that $\phi_1(n)$ will not change sign on the $S$-scale, but will have a sign change on the $l$-scale, as can be seen from the last factor in the right hand side of (\ref{s32_phil2}). For general $j\ge 1$, $\mathrm{Ai}(S+r_{j-1})$ will have $j-1$ zeros in the range $S>0$ and $\phi_j(n)$ will have an additional sign change on the $l$-scale, so that when $\rho>1$, as was the case when $\rho<1$, $\phi_j(n)$ will have exactly $j$ sign changes. The $O(K^{-1/3})$ correction term to (\ref{s32_phiS2}) is given in (\ref{s51_Phi(1)}) with $r_j$ replaced by $r_{j-1}$. We will again need different expansions on the $\xi$ and $n=O(1)$ scales, but for the former the $\phi_j(n)$ will be as in (\ref{s32_phixi}) with $j$ replaced by $j-1$, while for $n=O(1)$ (\ref{s32_phio1}) will hold, with ${k_1}$ and ${k_2}$ related by (\ref{s32_k1andk2}) with $r_j$ replaced by $r_{j-1}$.

\subsection{The case $\rho\sim 1$}\label{section:3.3}

Now we consider the case where $\rho\sim 1$, introducing the parameter $\eta$ by
\begin{equation}\label{s32_eta}
\rho=1+\frac{\eta}{K^{2/3}},\quad\textrm{i.e.,  }\eta=K^{2/3}(\rho-1)
\end{equation}
where $\eta$ can have either sign. Now we find that the eigenvalues behave as
\begin{equation}\label{s32_nueta}
\nu_j=\frac{1}{K}+\Big(\frac{\eta^2}{4}-r_j^*\Big)\frac{1}{K^{4/3}}+O(K^{-5/3}),\; j\ge 0
\end{equation}
where the $r_j^*$ are solutions to 
\begin{equation}\label{s32_rj*}
\mathrm{Ai}'(r_j^*)+\frac{\eta}{2}\mathrm{Ai}(r_j^*)=0.
\end{equation}
These roots depend upon $\eta$ so we write $r_j^*=r_j^*(\eta)$. We have $r_j^*<0$ for $j\ge 1$ while $r_0^*$ may have either sign. We order the roots as $r_0^*>r_1^*>r_2^*>\cdots$ so that $|r_j^*|<|r_{j+1}^*|$ for $j\ge 1$. Note that if $\eta=0$ (i.e., $\rho=1$) these are simply the roots of the derivative of the Airy function. 

The eigenvectors are concentrated on the scale $S=O(1)$ and we shall obtain
\begin{equation}\label{s32_phiS*}
\phi_j(n)\sim {k_0}\rho^{-n/2}\mathrm{Ai}(S+r_j^*(\eta)),\; j\ge 0.
\end{equation}
This applies both for $S=O(1)$ and for $l=K-n=O(1)$, since unlike (\ref{s32_phiS}), (\ref{s32_phiS*}) does not vanish as $S\to 0$. A correction term to (\ref{s32_phiS*}) is given in subsection \ref{section:4.3} (see (\ref{s51_Phiexpand}) and (\ref{s51_Phi(1)rho=1})). Again different expansions must be constructed on the $\xi$-scale, where
\begin{equation}\label{s32_phixi*}
\phi_j(n)\sim {k_1}\rho^{-n/2}K^{-1/12}L(\xi)e^{\sqrt{K}\psi(\xi)}e^{K^{1/6}\psi_j^*(\xi)},
\end{equation}
with $L(\cdot)$ and $\psi(\cdot)$ as in (\ref{s32_psixi}) and (\ref{s32_lxi}), and 
\begin{equation}\label{s32_psixi*}
\psi_j^*(\xi)=-\frac{1}{2}r_j^*(\eta)\Big[\sqrt{\xi(1-\xi)}+\sin^{-1}\big(\sqrt{1-\xi}\big)\Big].
\end{equation}
The constants ${k_0}$ and ${k_1}$ are related as in (\ref{s32_k0andk1}). For $n=O(1)$, (\ref{s32_phio1}) applies and  ${k_1}$ and ${k_2}$ are related as in (\ref{s32_k1andk2}), but with $r_j$ replaced by $r_j^*(\eta)$.  We also note that with the scaling in (\ref{s32_eta}) we can use
$\rho^{-n/2}\sim\exp\big(-\eta K^{1/3}/2+\eta S/2\big)$
to replace the $\rho^{-n/2}$ factors in (\ref{s32_phiS*}) and (\ref{s32_phixi*}). 

To better see the asymptotic matching between the results for $\rho-1=O(K^{-2/3})$ and those for $\rho\gtrless 1$, we analyze the solutions to (\ref{s32_rj*}) in the limits of $\eta\to\pm\infty$. For $\eta\to-\infty$ we have $r_j^*(-\infty)=r_j$ for $j\ge 0$, while for $\eta\to+\infty$ we have $r_j^*(+\infty)=r_{j-1}$ for $j\ge 1$. When $j=0$ we can show from (\ref{s32_rj*}) that 
$$\frac{\eta^2}{4}-r_0^*(\eta)\sim\frac{1}{\eta},\;\eta\to+\infty$$
so that the $j=0$ eigenvalue in (\ref{s32_nueta}) begins to resemble the expression in (\ref{s32_nu0}).

\subsection{Orthogonality relations}\label{section:3.4}

We can show, using (\ref{s2_fs_pn't})--(\ref{s2_fs_pn'tb2}), that $\phi_j(n)$ satisfy the (exact) orthogonality relation
\begin{equation}\label{s32_orth}
\sum_{n=0}^{K-1}\rho^n(n+1)\phi_i(n)\phi_j(n)=K_*(i;K,\rho)\delta_{ij},
\end{equation}
and $K_*$ is a normalizing constant (obtained by setting $i=j$ in (\ref{s32_orth})).
Using (\ref{s32_orth}) and (\ref{s2_fs_phi}), we obtain the explicit expression for the $c_j$ in (\ref{s2_cj}).

For $i,\,j\ge 0$ and $\rho<1$, and also for $i,\,j\ge 1$ and $\rho>1$, the relation in (\ref{s32_orth}) asymptotically reduces to 
\begin{equation*}\label{s32_orthairy}
\int_0^\infty \mathrm{Ai}(S+r_i)\mathrm{Ai}(S+r_j)dS=\big[\mathrm{Ai}'(r_j)\big]^2\delta_{ij}.
\end{equation*}
When $i=0$ and $j\ge 1$ with $\rho>1$ we can use (\ref{s32_phi0xi}) to approximate $\phi_0(n)$ and (\ref{s32_phil2}) to approximate $\phi_j(n)$ for $l=K-n=O(1)$. Since $\rho^n\phi_j(n)$ is concentrated in the range $n=K-O(1)$ for $\rho>1$, and $\phi_0(n)$ varies weakly with $n$ (as $\xi=n/K$) we can approximate $\phi_0(n)$ by a constant near $\xi=1$ and (\ref{s32_orth}) holds since
\begin{equation*}\label{s32_orth0}
\sum_{l=1}^\infty \rho^{-l/2}\Big(l-\frac{\sqrt{\rho}}{\sqrt{\rho}-1}\Big)=0,
\end{equation*}
which illustrates how $\phi_0$ is orthogonal to the other $\phi_j$ when $\rho>1$. When $\rho-1=\eta K^{-2/3}$, (\ref{s32_orth}) asymptotically reduces to 
\begin{equation*}\label{s32_orthreduce}
\int_0^\infty \mathrm{Ai}(S+r_i^*(\eta))\mathrm{Ai}(S+r_j^*(\eta))dS=\Big[\frac{\eta^2}{4}-r_i^*(\eta)\Big]\Big[\mathrm{Ai}(r_i^*(\eta))\Big]^2\delta_{ij}
\end{equation*}
and this relation can also be obtained directly using Sturm-Liouville theory for ODEs.

\subsection{Tail behaviors of the unconditional density}\label{section:3.5}

Finally, we give the asymptotic approximations to the unconditional sojourn time density $p(t)$ for sufficiently large $t$ (typically $t\gg O(K^{4/3})$). From this we shall see a variety of different tail behaviors for the finite capacity model.

For $\rho<1$, we have
\begin{eqnarray}\label{S35_rho<1}
p(t)&\sim& \frac{1+\sqrt{\rho}}{1-\sqrt{\rho}}\,e^{\frac{1+\sqrt{\rho}}{1-\sqrt{\rho}}}\big[\mathrm{Ai}'(r_0)\big]^{-2}\,K^{-4/3}\exp\Big\{-\pi\sqrt{K}-\frac{\pi}{2} r_0K^{1/6}\Big\}\nonumber\\
&&\times\exp\Big\{t\Big[-(1-\sqrt{\rho})^2-\frac{\sqrt{\rho}}{K}+\frac{\sqrt{\rho}}{K^{4/3}}r_0+O(K^{-5/3})\Big]\Big\},
\end{eqnarray}
where $r_0\approx -2.338$ is the largest root of $\mathrm{Ai}(r)=0$.

When $\rho\sim 1$, we need to consider four different scales, which are $\rho=1+O(K^{-1/2})<1$, $\rho=1+O(K^{-2/3})<1$, $\rho=1+O(K^{-1})$, and $\rho=1+O(K^{-2/3})>1$.

We first introduce the parameter $b$ by setting $\rho=1+b/\sqrt{K}$. For $b<0$, we have the asymptotic approximation
\begin{eqnarray}\label{S35_rho1_b}
p(t) &\sim & K^{-5/6}\frac{4|b|}{b^2+4}\big[\mathrm{Ai}'(r_0)\big]^{-2}\exp\bigg\{\frac{2}{b}\sin^{-1}\Big(\frac{|b|}{\sqrt{b^2+4}}\Big)-\frac{2b^2+4}{b^2+4}\bigg\}\nonumber\\
&&\times\exp\bigg\{-2\sin^{-1}\Big(\frac{|b|}{\sqrt{b^2+4}}\Big)\sqrt{K}+K^{1/6}r_0\Big[\frac{2b}{b^2+4}-\sin^{-1}\Big(\frac{|b|}{\sqrt{b^2+4}}\Big)\Big]\bigg\}\nonumber\\
&&\times\exp\bigg\{t\Big[-\Big(1+\frac{b^2}{4}\Big)\frac{1}{K}+\frac{r_0}{K^{4/3}}+O(K^{-5/3})\Big]\bigg\}.
\end{eqnarray}
This holds for $K^{4/3}\ll t\ll K^{5/3}$, and can be extended to even larger time ranges by using the higher terms in (\ref{s32_nu}).

Then we consider the $\eta$-scale as in (\ref{s32_eta}) with $\eta=K^{2/3}(\rho-1)<0$, and the asymptotic approximation is  
\begin{eqnarray}\label{S35_rho1_eta-}
p(t) &\sim& \frac{|\eta|}{K}\frac{\Big(\int_0^\infty e^{-\eta S/2}\mathrm{Ai}(S+r_0^*(\eta))\,dS\Big)^2}{\Big({\eta^2}/{4}-r_0^*(\eta)\Big)\big[\mathrm{Ai}(r_0^*(\eta))\big]^2}\,\exp\bigg\{\eta K^{1/3}\bigg\}\nonumber\\
&&\times\exp\bigg\{t\Big[-\frac{1}{K}-\frac{{\eta^2}/{4}-r_0^*(\eta)}{K^{4/3}}+O(K^{-5/3})\Big]\bigg\},
\end{eqnarray}
where $r_0^*(\eta)$ is the solution to (\ref{s32_rj*}) as $j=0$.

Next we introduce the parameter $a$ by $\rho=1+a/K$, with $-\infty<a<\infty$. With this scaling we have the following approximation
\begin{equation}\label{S35_rho1_a}
p(t) \sim \frac{a\,e^a}{e^a-1}\frac{\Big(\int_0^\infty \mathrm{Ai}(S+r_0^*(0))\,dS\Big)^2}{|r_0^*(0)|\big[\mathrm{Ai}(r_0^*(0))\big]^2}K^{-4/3}\exp\bigg\{t\Big[-\frac{1}{K}-\frac{|r_0^*(0)|}{K^{4/3}}+O(K^{-5/3})\Big]\bigg\},
\end{equation}
where $r_0^*(0)\approx -1.019$ is the largest root of $\mathrm{Ai}'(r)=0$.

Under the $\eta$-scale in (\ref{s32_eta}) with $\eta>0$, we have 
\begin{equation}\label{S35_rho1_eta+}
p(t) \sim \frac{\eta}{K}\frac{\Big(\int_0^\infty e^{-\eta S/2}\mathrm{Ai}(S+r_0^*(\eta))\,dS\Big)^2}{\Big(\eta^2/4-r_0^*(\eta)\Big)\big[\mathrm{Ai}(r_0^*(\eta))\big]^2}\exp\bigg\{t\Big[-\frac{1}{K}-\frac{{\eta^2}/{4}-r_0^*(\eta)}{K^{4/3}}+O(K^{-5/3})\Big]\bigg\}.
\end{equation}

For $\rho>1$, the tail of $p(t)$ behaves as the exponential density on the time scale $t=O(K)$, with
\begin{equation}\label{S35_rho>1}
p(t)\sim \frac{1}{K}\,e^{-t/K}.
\end{equation}
Now (\ref{S35_rho>1}) holds for all $t=O(K)$ and most of the probability mass concentrates on this time scale. Note that by using the correction terms in (\ref{s32_nu}) we can extend the validity of (\ref{S35_rho>1}) to larger time scales $t\gg O(K)$.

From the approximations in (\ref{S35_rho1_eta-})--(\ref{S35_rho>1}), we can form the uniform asymptotic approximation 
\begin{eqnarray}\label{S35_uniform}
p(t) &\sim& \frac{\eta}{K}\frac{\Big(\int_0^\infty e^{-\eta S/2}\mathrm{Ai}(S+r_0^*(\eta))\,dS\Big)^2}{\Big({\eta^2}/{4}-r_0^*(\eta)\Big)\big[\mathrm{Ai}(r_0^*(\eta))\big]^2}\frac{\exp\Big(\eta K^{1/3}\Big)}{\exp\Big(\eta K^{1/3}\Big)-1}\nonumber\\
&&\times\exp\bigg\{t\Big[-\frac{1}{K}-\frac{{\eta^2}/{4}-r_0^*(\eta)}{K^{4/3}}+O(K^{-5/3})\Big]\bigg\},
\end{eqnarray}
which holds for any $-\infty<\eta\le \infty$ and any $a=\eta K^{1/3}$. The results in (\ref{S35_rho1_eta-})--(\ref{S35_rho>1}) are limiting cases of (\ref{S35_uniform}). But we still need to use the result in (\ref{S35_rho<1}) for $\rho<1$, and in (\ref{S35_rho1_b}) for $\rho\sim 1$ and $-\infty<b<0$ (corresponding to large negative $\eta$ values).

\section{Brief derivations}\label{section:4}

We analyze the spectrum of the finite capacity model for $K\to\infty$. In \cite{ZHE_O09} we give an explicit expression for the Laplace transform, $\int_0^\infty e^{-\theta t}p_n(t)dt$, of the conditional sojourn time density. But this involves a complicated expression with integrals related to hypergeometric functions. The poles of this expression in the complex $\theta$-plane correspond to the (exact) eigenvalues, while the residues at the poles yield the eigenvectors. We did compute the dominant eigenvalue $\nu_0$ in \cite{ZHE_O09}, for $K\to\infty$ and the three cases of $\rho$. However, the eigenvalues $\nu_j$ for $j\ge 1$ were not obtained, nor the eigenvectors $\phi_j(n)$ for any $j\ge 0$. While these could in principle be obtained by expanding the exact expression in \cite{ZHE_O09}, here we instead apply a singular perturbation approach. We shall thus obtain from (\ref{s2_fs_pn't})-(\ref{s2_fs_pkt}) limiting differential equations with appropriate boundary conditions. This analysis is much simpler than trying to obtain the results from \cite{ZHE_O09}, and is also applicable to problems for which no exact solution is available.

\subsection{The case $\rho<1$}\label{section:4.1}
We first consider first the scale $n=K-O(K^{2/3})$ and introduce the new variable $S$, with 
\begin{equation}\label{s51_Sscale}
n=K-K^{2/3}S,\quad\textrm{i.e. } S=\frac{K-n}{K^{2/3}}\in(0,\infty).
\end{equation} 
We shall see that the analysis ceases to be valid for $S\to 0$ and also $S\to\infty$, and later we analyze separately the scales $n=K-O(1)$ (then $S=O(K^{-2/3})$) and $n=O(K)$ (then $S=O(K^{1/3})$). But, the scale $S=O(1)$ is the one which governs the sign changes of the eigenvectors and leads to a limiting ODE.

Setting $p_n(t)=e^{-\nu t}\rho^{-n/2}\Phi(S)$ we obtain from (\ref{s2_fs_pn't})
\begin{equation}\label{s51_Phirec}
-\nu\Phi(S)=\sqrt{\rho}\big[\Phi(S-K^{-2/3})+\Phi(S+K^{-2/3})\big]-(1+\rho)\Phi(S)-\frac{\sqrt{\rho}}{K-K^{2/3}S+1}\Phi(S+K^{-2/3}),
\end{equation}
where we note that changing $n\to n\pm 1$ corresponds to, in view of (\ref{s51_Sscale}), changing $S\to S\mp K^{-2/3}$. We use the artificial boundary condition in (\ref{s2_fs_pkt}), which implies that
\begin{equation}\label{s51_Phi(0)}
\Phi(0)=\sqrt{\rho}\,\Phi(K^{-2/3}).
\end{equation}
From (\ref{s51_Phirec}), letting $K\to\infty$ leads to the conclusion that $\nu\to 1+\rho-2\sqrt{\rho}=(1-\sqrt{\rho})^2$, which is the relaxation rate for the standard $M/M/1$ queue, so we see a coalescence of the eigenvalues. Also, from (\ref{s51_Phi(0)}) we conclude that, if $\rho\ne 1$, $\Phi(0)=0$, at least to leading order for $K$ large.

To proceed further we write
\begin{equation}\label{s51_nuexpand}
\nu=(1-\sqrt{\rho})^2+\frac{\sqrt{\rho}}{K}+\frac{\sqrt{\rho}}{K^{4/3}}\widetilde{\nu}
\end{equation}
and then expand the eigenvalues and eigenfunctions as
\begin{equation*}
\widetilde{\nu}=\widetilde{\nu}^{(0)}+\frac{\widetilde{\nu}^{(1)}}{K^{1/3}}+\frac{\widetilde{\nu}^{(2)}}{K^{2/3}}+O(K^{-1}),
\end{equation*}
and
\begin{equation}\label{s51_Phiexpand}
\Phi(S)=\Phi^{(0)}(S)+\frac{1}{K^{1/3}}\Phi^{(1)}(S)+O(K^{-2/3}).
\end{equation}
We do not indicate above the dependence on the eigenvalue index $j$. The expansions in (\ref{s51_nuexpand}) and (\ref{s51_Phiexpand}) are necessary to obtain from (\ref{s51_Phirec}) and (\ref{s51_Phi(0)}) a limiting ODE and boundary condition. Note that (\ref{s51_Phirec}) can be used to conclude that $\nu-(1-\sqrt{\rho})^2\sim\sqrt{\rho}K^{-1}$ as $K\to\infty$, so that the first two terms in (\ref{s51_nuexpand}) are independent of $j$. 

Using (\ref{s51_nuexpand}) and (\ref{s51_Phiexpand}) in (\ref{s51_Phirec}), we obtain at leading order ($O(K^{-4/3}))$
\begin{equation}\label{s51_Phi0ODE}
\frac{d^2}{dS^2}\Phi^{(0)}(S)-(S-\widetilde{\nu}^{(0)})\Phi^{(0)}(S)=0
\end{equation}
and at the next order ($O(K^{-5/3})$) we get
\begin{equation}\label{s51_Phi1ODE}
\frac{d^2}{dS^2}\Phi^{(1)}(S)-(S-\widetilde{\nu}^{(0)})\Phi^{(1)}(S)=-\widetilde{\nu}^{(1)}\Phi^{(0)}(S)+S^2\Phi^{(0)}(S)+\frac{d}{dS}\Phi^{(0)}(S).
\end{equation}
From (\ref{s51_Phi(0)}) we obtain the boundary conditions 
\begin{equation}\label{s51_ICphi0phi1}
\Phi^{(0)}(0)=0,\quad \Phi^{(1)}(0)=0.
\end{equation}
We also require that $\Phi^{(l)}(S)$ decay as $S\to\infty$. Equation (\ref{s51_Phi0ODE}) is the Airy equation and the decaying solution is 
$\Phi^{(0)}(S)=\mathrm{Ai}(S-\widetilde{\nu}^{(0)}).$
But then $\Phi^{(0)}(0)=0$ implies that $\mathrm{Ai}(-\widetilde{\nu}^{(0)})=0$ so that $-\widetilde{\nu}^{(0)}$ must be a root of the Airy function. We denote the roots of $\mathrm{Ai}(\cdot)$ as $r_0>r_1>r_2>\cdots$ and we have $r_j<0$. We have thus obtained $\widetilde{\nu}^{(0)}=-r_j=|r_j|$ and, up to a constant, $\Phi^{(0)}(S)=\mathrm{Ai}(S+r_j)$, and thus derived (\ref{s32_phiS}) and the $O(K^{-4/3})$ term(s) in (\ref{s32_nu}). 

We can easily compute higher order terms in the expansions in (\ref{s51_nuexpand}) and (\ref{s51_Phiexpand}). For example, $\widetilde{\nu}^{(1)}$, which gives the $O(K^{-5/3})$ term(s) in the expansion(s) of the eigenvalues, can be obtained from the solvability condition for (\ref{s51_Phi1ODE}). We multiply (\ref{s51_Phi1ODE}) by $\Phi^{(0)}(S)$ and integrate over $S\in (0,\infty)$. Then some integration by parts and use of (\ref{s51_ICphi0phi1}) leads to 
\begin{equation}\label{s51_solvab}
\widetilde{\nu}^{(1)}\int_0^\infty\big[\Phi^{(0)}(S)\big]^2dS=\int_0^\infty S^2\big[\Phi^{(0)}(S)\big]^2dS.
\end{equation} 
Using results in \cite{ALB}, namely, 
\begin{equation}\label{s51_alb01}
\int \big[\mathrm{Ai}(S)\big]^2dS=S\big[\mathrm{Ai}(S)\big]^2-\big[\mathrm{Ai}'(S)\big]^2,
\end{equation}
\begin{equation}\label{s51_alb02}
\int S\big[\mathrm{Ai}(S)\big]^2dS=\frac{1}{3}\Big\{\mathrm{Ai}(S)\mathrm{Ai}'(S)-S\big[\mathrm{Ai}'(S)\big]^2+S^2\big[\mathrm{Ai}(S)\big]^2\Big\},
\end{equation}
\begin{equation}\label{s51_alb03}
\int S^2\big[\mathrm{Ai}(S)\big]^2dS=\frac{1}{5}\Big\{2S\mathrm{Ai}(S)\mathrm{Ai}'(S)-\big[\mathrm{Ai}(S)\big]^2-S^2\big[\mathrm{Ai}'(S)\big]^2+S^3\big[\mathrm{Ai}(S)\big]^2\Big\},
\end{equation}
and the fact that $\mathrm{Ai}(r_j)=0$, we find that 
\begin{equation*}\label{s51_airy01}
\int_0^\infty \big[\mathrm{Ai}(S+r_j)\big]^2dS=\big[\mathrm{Ai}'(r_j)\big]^2,
\end{equation*}
\begin{equation*}\label{s51_airy02}
\int_0^\infty S^2\big[\mathrm{Ai}(S+r_j)\big]^2dS=\frac{8}{15}r_j^2\big[\mathrm{Ai}'(r_j)\big]^2,
\end{equation*}
which yield $\widetilde{\nu}^{(1)}=8r_j^2/15$, and this establishes the $O(K^{-5/3})$ term in (\ref{s51_nuexpand}) and (\ref{s32_nu}). Given $\widetilde{\nu}^{(1)}$ the correction term $\Phi^{(1)}$ in (\ref{s51_Phiexpand}) can be obtained by solving the inhomogeneous Airy equation in (\ref{s51_Phi1ODE}). We can construct a particular solution to (\ref{s51_Phi1ODE}) in the form $\Phi^{(1)}=u_1(S)\mathrm{Ai}(S+r_j)+u_2(S)\mathrm{Bi}(S+r_j)$, where
\begin{equation*}
u_1(S)=-\pi\int\Big(\mathrm{Ai}'(S+r_j)\mathrm{Bi}(S+r_j)+(S^2-\widetilde{\nu}^{(1)})\,\mathrm{Ai}(S+r_j)\mathrm{Bi}(S+r_j)\Big)\,dS,
\end{equation*}
\begin{equation*}
u_2(S)=\pi\int\Big(\mathrm{Ai}(S+r_j)\mathrm{Ai}'(S+r_j)+(S^2-\widetilde{\nu}^{(1)})\,\mathrm{Ai}^2(S+r_j)\Big)\,dS.
\end{equation*}
Again, by using results in \cite{ALB} and some further calculations, we have
\begin{equation}\label{s51_Phi(1)}
\Phi^{(1)}(S)=\Big(\frac{3}{10}S+\frac{19}{30}r_j\Big)\mathrm{Ai}(S+r_j)+\Big(\frac{1}{5}S^2-\frac{4r_j}{15}S\Big)\mathrm{Ai}'(S+r_j).
\end{equation}
Note that to this particular solution we can add an arbitrary multiple of the homogeneous solution, say $k_0^{(1)}\mathrm{Ai}(S+r_j)$. But if we allow $k_0$ in (\ref{s32_phiS}) to depend upon $K$, say via an expansion in powers of $K^{-1/3}$, then we can set $k_0^{(1)}=0$ and view the homogeneous solution as corresponding to the correction term in the expansion of $k_0$. This completes the analysis of the $S$-scale.

Now, the approximation $\rho^{n/2}\phi_j(n)\sim {k}_0\mathrm{Ai}(S+r_j)$ has $j$ zeros in the range $S>0$, at $S=r_l-r_j=|r_l-r_j|$ for $l=0,1,\cdots,j-1$. This corresponds to $j$ sign changes in the $\phi_j(n)$, which is to be expected. However, as $S\to 0$, $\mathrm{Ai}(S+r_j)\to 0$ as we must examine carefully the range where $S$ is small, in particular to see if other sign changes occur. We consider $n=K-O(1)$ so we set $l=K-n$. Then on the $l=O(1)$ scale we let
\begin{equation}\label{s51_phij(n)}
\phi_j(n)=\rho^{-n/2}{k}_1K^{-2/3}\mathcal{L}(l;K)
\end{equation}  
where from (\ref{s2_fs_pn't}) we find that $\mathcal{L}$ satisfies
\begin{equation*}\label{s51_mathcalS}
-\nu\mathcal{L}(l)=\sqrt{\rho}\big[\mathcal{L}(l-1)+\mathcal{L}(l+1)\big]-(1+\rho)\mathcal{L}(l)-\frac{\sqrt{\rho}}{K}\Big(1-\frac{l-1}{K}\Big)^{-1}\mathcal{L}(l)
\end{equation*}
and (\ref{s2_fs_pn'tb2}) leads to $\mathcal{L}(0)=\sqrt{\rho}\mathcal{L}(1)$. Given $\nu\sim (1-\sqrt{\rho})^2$ we see to leading order that $2\mathcal{L}^{(0)}(l)=\mathcal{L}^{(0)}(l-1)+\mathcal{L}^{(0)}(l+1)$, where $\mathcal{L}^{(0)}$ is the leading term in the expansion of $\mathcal{L}$ as $K\to\infty$. Then also $\mathcal{L}^{(0)}(0)=\sqrt{\rho}\mathcal{L}^{(0)}(1)$ and hence, up to a constant that can be incorporated into ${k}_1$ in (\ref{s51_phij(n)}), we have 
\begin{equation}\label{s51_s(0)l}
\mathcal{L}^{(0)}(l)=l+\frac{\sqrt{\rho}}{1-\sqrt{\rho}};\quad l=0,1,\cdots.
\end{equation}
By asymptotic matching the behavior of the right side of (\ref{s51_phij(n)}) as $l\to\infty$ must agree with the expansion of (\ref{s32_phiS}) as $S\to 0$, and thus ${k}_1K^{-2/3}l\sim{k}_0\mathrm{Ai}'(r_j)S={k}_0\mathrm{Ai}'(r_j)K^{-2/3}l$ so that ${k}_1={k}_0\mathrm{Ai}'(r_j)$, and we have derived (\ref{s32_phil}). Note that (\ref{s51_s(0)l}) is strictly positive for $l\ge 0$ if $\rho<1$. Thus all of the sign changes in the eigenvectors $\phi_j(n)$ occur on the $S$-scale if $\rho<1$. 

We next consider the scales $n=O(K)$ and $n=O(1)$, thus obtaining the ``tails'' of the eigenfunctions. On the former scale we set
\begin{equation}\label{s51_xiscale}
\xi=\frac{n}{K},\quad \phi_j(n)=\rho^{-n/2}\widetilde{\Phi}(\xi;K,\rho)
\end{equation}
with which (\ref{s2_fs_pn't}) becomes
\begin{equation}\label{s51_Phitilderec}
-\nu\widetilde{\Phi}(\xi)=\sqrt{\rho}\Big[\widetilde{\Phi}\Big(\xi+\frac{1}{K}\Big)+\widetilde{\Phi}\Big(\xi-\frac{1}{K}\Big)\Big]-(1+\rho)\widetilde{\Phi}(\xi)-\frac{\sqrt{\rho}}{K\xi}\Big(1+\frac{1}{K\xi}\Big)^{-1}\widetilde{\Phi}\Big(\xi-\frac{1}{K}\Big),
\end{equation}
where we again do not indicate the dependence of $\widetilde{\Phi}$ on the eigenvalue index $j$. 

If we examine the result on the $S$-scale for large $S$ we have, apart from the $\rho^{-n/2}$ factor, 
\begin{eqnarray}\label{s51_largeS}
{k}_0\mathrm{Ai}(S+r_j)&\sim& \frac{{k}_0}{2\sqrt{\pi}}S^{-1/4}\exp\Big(-\frac{2}{3}S^{3/2}-r_j\sqrt{S}\Big)\nonumber\\
&=&\frac{{k}_0}{2\sqrt{\pi}}K^{-1/12}(1-\xi)^{-1/4}\exp\Big(-\frac{2}{3}\sqrt{K}(1-\xi)^{3/2}-r_jK^{1/6}\sqrt{1-\xi}\Big),
\end{eqnarray}
where we used $S=K^{1/3}(1-\xi)$. If the $\xi$- and $S$-scale results will match asymptotically, the behavior of $\widetilde{\Phi}$ as $\xi\uparrow 1$ must be of the form in the right side of (\ref{s51_largeS}). Then we assume the WKB-type ansatz
\begin{equation}\label{s51_WKB}
\widetilde{\Phi}(\xi;K,\rho)={k}_1K^{-1/12}\exp\Big\{\sqrt{K}\psi(\xi)+K^{1/6}\psi^{(1)}(\xi)\Big\}\Big[L(\xi)+K^{-1/6}L^{(1)}(\xi)+\cdots\Big].
\end{equation}
The constants ${k}_0$ and ${k}_1$ will be related by (\ref{s32_k0andk1}), in view of (\ref{s51_largeS}) and the matching condition. 

Using (\ref{s51_WKB}) in (\ref{s51_Phitilderec}) we obtain at the first three orders ($O(K^{-1})$, $O(K^{-4/3})$ and $O(K^{-3/2})$) the ODEs
\begin{equation}\label{s51_order1}
\big[\psi'(\xi)\big]^2=\frac{1}{\xi}-1,
\end{equation}
\begin{equation}\label{s51_order2}
2(\psi^{(1)}(\xi))'\psi'(\xi)=-r_j,
\end{equation}
\begin{equation}\label{s51_order3}
2\psi'(\xi)L'(\xi)+\Big[\psi''(\xi)+\frac{1}{\xi}\psi'(\xi)\Big]L(\xi)=0.
\end{equation}
From (\ref{s51_largeS}) we also see that $\psi(1)=\psi^{(1)}(1)=0$ and more precisely, as $\xi\uparrow 1$, $\psi(\xi)\sim-\frac{2}{3}(1-\xi)^{3/2}$ and $\psi^{(1)}(\xi)\sim -r_j\sqrt{1-\xi}$. Then integrating (\ref{s51_order1}) yields
\begin{equation}\label{s51_phiint}
\psi(\xi)=-\int_\xi^1\sqrt{\frac{1-v}{v}}dv=\sqrt{\xi(1-\xi)}-\sin^{-1}(\sqrt{1-\xi})
\end{equation}
with $\sin^{-1}(\cdot)\in[0,\pi/2]$. Given (\ref{s51_phiint}) we can easily integrate (\ref{s51_order2}) and (\ref{s51_order3}) to get (\ref{s32_psixi1}) and (\ref{s32_lxi}).

Finally we consider the scale $n=O(1)$. This is necessary since $L(\xi)$ in (\ref{s32_lxi}) is singular as $\xi\to 0$, behaving as $L(\xi)\sim \xi^{-1/4}$. Setting $1-\xi=1-n/K$ and expanding (\ref{s32_phixi}) (or (\ref{s51_WKB})) for $\xi\to 0$ leads to 
\begin{equation}\label{s51_expand360}
\rho^{-n/2}{k}_1K^{1/6}n^{-1/4}\exp\Big\{\sqrt{K}\psi(0)+2\sqrt{n}+K^{1/6}\psi^{(1)}(0)\Big\},
\end{equation}
where we used $\sin^{-1}(\sqrt{1-\xi})\sim \pi/2-\sqrt{\xi},\;\xi\to 0$, and note that $\psi(0)=-\pi/2$ and $\psi^{(1)}(0)=-\pi r_j/4$. For $n=O(1)$ we set $\phi_j(n)\sim \rho^{-n/2}Q(n)$ where, from (\ref{s2_fs_pn't}), we find that 
$$Q(n+1)+Q(n-1)-2Q(n)=\frac{1}{n+1}Q(n-1),$$
whose solution is the contour integral in (\ref{s32_phio1}). Then equating the large $n$ behavior of this integral to (\ref{s51_expand360}) gives, by matching, the relation in (\ref{s32_k1andk2}). This completes the analysis of the case $\rho<1$.

\subsection{The case $\rho>1$}\label{section:4.2}
When $\rho>1$ and $K\to\infty$ we expect that the number of customers in the PS queue will typically be close to the capacity $K$. We then note that the analysis of the case $\rho<1$ made use of the fact that $\rho\ne 1$ in order to conclude from (\ref{s51_Phi(0)}) that $\Phi(0)=0$. Note also that if $\rho=1$, (\ref{s51_Phi(0)}) would imply, to leading order, that $\Phi'(0)=0$. The calculations in subsection \ref{section:4.1} apply equally well to the case $\rho>1$, for each of the four ranges of $n$. The one crucial difference, however, is that if $\rho>1$ the expansion of the $\phi_j(n)$ on the scale $l=K-n=O(1)$, cf. (\ref{s51_s(0)l}), undergoes a sign change as $l$ increases from $l=1$ to $l=\infty$. But $\phi_j(n)$ should have exactly $j$ sign changes with $n$. We thus conclude that when $\rho>1$ the expansion of $\phi_j(n)$ in subsection \ref{section:4.1} corresponds to the $(j+1)^{st}$ eigenvector, and the expansion of $\nu_j$ to the $(j+1)^{st}$ eigenvalue, for $j\ge 0$. This leads to (\ref{s32_phiS2}) and (\ref{s32_phil2}). 

It remains to compute $\phi_0(n)$ and $\nu_0$, for $K\to\infty$ and $\rho>1$. We recall that for the finite population model \cite{ZHE_O12}, when the corresponding $\rho$ exceeded unity, we saw all of the eigenfunctions varying smoothly on the $\xi$-scale, but their zeros concentrated near $\xi=1-\rho^{-1}$. We thus consider here the scale $\xi=n/K\in (0,1)$, setting $p_n(t)=e^{-\nu t}\varphi(\xi)$, and then (\ref{s2_fs_pn't}) becomes
\begin{equation}\label{s52_phirec}
-\nu\varphi(\xi)=\rho\Big[\varphi\Big(\xi+\frac{1}{K}\Big)-\varphi(\xi)\Big]+\varphi\Big(\xi-\frac{1}{K}\Big)-\varphi(\xi)-\frac{1}{K\xi}\Big(1+\frac{1}{K\xi}\Big)^{-1}\varphi\Big(\xi-\frac{1}{K}\Big).
\end{equation}
Note that (\ref{s52_phirec}) differs from (\ref{s51_Phitilderec}) as the present analysis does not involve the symmetrizing factor $\rho^{-n/2}$. Since the right side of (\ref{s52_phirec}) is approximately $\big[(\rho-1)\varphi'(\xi)-\xi^{-1}\varphi(\xi)\big]K^{-1}+O(K^{-2})$, we expand the eigenvalues and eigenfunctions  as
\begin{equation*}\label{s52_nu}
\nu=\frac{1}{K}\nu^{(1)}+\frac{1}{K^2}\nu^{(2)}+O(K^{-3}),
\end{equation*}
\begin{equation}\label{s52_phi}
\varphi(\xi)=\varphi^{(1)}(\xi)+\frac{1}{K}\varphi^{(2)}(\xi)+O(K^{-2}).
\end{equation} 
Then (\ref{s52_phirec}) leads to, at orders $O(K^{-1})$ and $O(K^{-2})$, 
\begin{equation}\label{s52_rec1}
-\nu^{(1)}\varphi^{(1)}(\xi)=(\rho-1)\frac{d}{d\xi}\varphi^{(1)}(\xi)-\frac{1}{\xi}\varphi^{(1)}(\xi)
\end{equation}
and 
\begin{eqnarray}\label{s52_rec2}
-\nu^{(2)}\varphi^{(1)}(\xi)-\nu^{(1)}\varphi^{(2)}(\xi)&=&(\rho-1)\frac{d}{d\xi}\varphi^{(2)}(\xi)-\frac{1}{\xi}\varphi^{(2)}(\xi)\nonumber\\
&&+\frac{1}{2}(\rho+1)\frac{d^2}{d\xi^2}\varphi^{(1)}(\xi)+\frac{1}{\xi}\frac{d}{d\xi}\varphi^{(1)}(\xi)+\frac{1}{\xi^2}\varphi^{(1)}(\xi).
\end{eqnarray}
The boundary condition in (\ref{s2_fs_pn'tb2}) leads to $\varphi(1)=\varphi(1-K^{-1})$ and with the expansion in (\ref{s52_phi}) we obtain
\begin{equation}\label{s52_phi1ODE}
\frac{d}{d\xi}\varphi^{(1)}(\xi)\Big|_{\xi=1}=0
\end{equation}
and 
\begin{equation}\label{s52_phi2ODE}
\Big[-\frac{d}{d\xi}\varphi^{(2)}(\xi)+\frac{1}{2}\frac{d^2}{d\xi^2}\varphi^{(1)}(\xi)\Big]\Big|_{\xi=1}=0.
\end{equation}

Solving the simple differential equation in (\ref{s52_rec1}) leads to, up to a multiplicative constant, 
\begin{equation}\label{s52_phi(1)}
\varphi^{(1)}(\xi)=\xi^{1/(\rho-1)}\exp\Big(-\frac{\nu^{(1)}}{\rho-1}\xi\Big)
\end{equation}
and then (\ref{s52_phi1ODE}) leads to $\nu^{(1)}=1$. This yields the leading term for the zeroth eigenvalue, and (\ref{s52_phi(1)}) shows that the eigenvector has no sign changes for $\xi\in(0,1)$. We then set
\begin{equation}\label{s52_phi(2)}
\varphi^{(2)}(\xi)=\xi^{1/(\rho-1)}\exp\Big(-\frac{\xi}{\rho-1}\Big)\overline{\varphi}(\xi)
\end{equation}
in (\ref{s52_rec2}) to obtain
\begin{equation}\label{s52_nu2}
-\nu^{(2)}=(\rho-1)\overline{\varphi}'(\xi)+\frac{1}{2}(\rho+1)\bigg[\frac{1}{(\rho-1)^2}\Big(\frac{1}{\xi}-1\Big)^2-\frac{1}{\rho-1}\frac{1}{\xi^2}\bigg]+\frac{1}{\rho-1}\Big(\frac{1}{\xi}-1\Big)\frac{1}{\xi}+\frac{1}{\xi^2}.
\end{equation}
From (\ref{s52_phi(1)}) and the fact $\nu^{(1)}=1$ we obtain $\frac{d^2}{d\xi^2}\varphi^{(1)}(\xi)\big|_{\xi=1}=-(\rho-1)^{-1}\varphi^{(1)}(1)$ and then (\ref{s52_phi2ODE}) and (\ref{s52_phi(2)}) show that 
\begin{equation}\label{s52_phibar'}
\overline{\varphi}'(1)=-\frac{1}{2(\rho-1)}.
\end{equation}
Setting $\xi=1$ in (\ref{s52_nu2}) and using (\ref{s52_phibar'}) yields $\nu^{(2)}=(\rho-1)^{-1}$, which corresponds to the $O(K^{-2})$ term in (\ref{s32_nu0}). In (\ref{s32_nu0}) we also gave the $O(K^{-3})$ term which we do not derive here; this would follow by examining the problem for $\varphi^{(3)}(\xi)$ and $\nu^{(3)}$ in (\ref{s52_phi}). Given $\nu^{(2)}$ the solution to (\ref{s52_nu2}) is
\begin{equation}\label{s52_phibar}
\overline{\varphi}(\xi)=-\frac{3\rho-1}{2(\rho-1)^3}\xi+\frac{\rho^2-\rho+2}{2(\rho-1)^3\xi}+\frac{2\rho}{(\rho-1)^3}\log\xi,
\end{equation}
and (\ref{s52_phi(2)}), (\ref{s52_phi}) and (\ref{s52_phibar}) give the $O(K^{-1})$ correction to the zeroth eigenvector $\phi_0(n)$ in (\ref{s32_phi0xi}).

When $\rho>1$, the leading term in (\ref{s52_phi(1)}) vanishes as $\xi\to 0^+$. This shows that the scale $n=O(1)$ must also be considered. We also note that if $\rho<1$ the expression in (\ref{s52_phi(1)}) has a singularity at $\xi=0$, and hence the construction of $\phi_0(n)$ done here can only be done for $\rho>1$. For $n=O(1)$ we find that $\phi_0(n)\sim q_n^{(0)}$ where 
$$0=\rho\Big(q_{n+1}^{(0)}-q_n^{(0)}\Big)+\frac{n}{n+1}\,q_{n-1}^{(0)}-q_n^{(0)}$$
with $q_{-1}^{(0)}$ finite. Solving the above difference equation using generating functions or contour integrals, we obtain the formula in (\ref{s32_phi0n}) as the approximation to $\phi_0(n)$ for $n=O(1)$. Then the relation in (\ref{s32_k0*andk1*}) follows from asymptotic matching of the $\xi$ and $n$ scale results. This completes the analysis for $\rho>1$.

\subsection{The case $\rho\sim 1$}\label{section:4.3}
When $\rho$ is close to unity, we must re-examine the boundary condition in (\ref{s2_fs_pkt}), which has the form (\ref{s51_Phi(0)}) on the $S$-scale. From (\ref{s51_Phi(0)}) we see that 
\begin{equation}\label{s53_Phirec}
0=(\sqrt{\rho}-1)\Phi(0)+\sqrt{\rho}K^{-2/3}\Phi'(0)+O(K^{-4/3}),
\end{equation}
so the first two terms balance when $\rho=1+O(K^{-2/3})$. We thus consider the scaling in (\ref{s32_eta}), and we will then consider the four spatial scales $K-n=K^{2/3}S=O(K^{2/3})$, $n=O(K)$ and $n=O(1)$. 

When $S=O(1)$ we again use the expansion in (\ref{s51_nuexpand}) and (\ref{s51_Phiexpand}) and obtain for $\Phi^{(0)}$ the Airy equation in (\ref{s51_Phi0ODE}). But, for $\rho-1=K^{-2/3}\eta$ we have $\sqrt{\rho}-1\sim\frac{1}{2}\eta K^{-2/3}$ and (\ref{s53_Phirec}) leads to 
\begin{equation}\label{s53_Phi0ODE}
\frac{d}{dS}\Phi^{(0)}(S)\Big|_{S=0}+\frac{\eta}{2}\Phi^{(0)}(0)=0,
\end{equation}
which is the boundary condition for the leading order eigenfunction approximation. We also note that if $\widetilde{\nu}\sim\widetilde{\nu}^{(0)}$, then, with the present scaling, (\ref{s51_nuexpand}) becomes
\begin{equation}\label{s53_nu}
\nu=\frac{1}{K}+\frac{1}{K^{4/3}}\Big(\frac{\eta^2}{4}+\widetilde{\nu}^{(0)}\Big)+O(K^{-5/3}).
\end{equation}
Solving (\ref{s51_Phi0ODE}) subject to (\ref{s53_Phi0ODE}) yields, up to a constant, 
${\Phi}^{(0)}(S)=\mathrm{Ai}(S+r_j^*),$
where $r_j^*=-\widetilde{\nu}^{(0)}$ are solutions to 
\begin{equation}\label{s53_rj*}
\mathrm{Ai}'(r_j^*)+\frac{\eta}{2}\mathrm{Ai}(r_j^*)=0.
\end{equation}
Thus $r_j^*=r_j^*(\eta)$ and if $\eta=0\;(\rho=1)$ these are the roots of $\mathrm{Ai}'(z)=0$. In Figure \ref{figure:r*} we include a sketch of the solution branches $r_j^*(\eta)$ of (\ref{s53_rj*}). We order the roots again as $r^*_0>r^*_1>r^*_2>\cdots$ and we have $r^*_j<0$ for $j\ge 1$. The zeroth root $r^*_0$ is negative for $\eta<-2\mathrm{Ai}'(0)/\mathrm{Ai}(0)=\pi^{-1}3^{5/6}\Gamma^2(2/3)=1.458\cdots$, but positive for $\eta>-2\mathrm{Ai}'(0)/\mathrm{Ai}(0)$. The roots satisfy the bounds 
$r_{j-1}>r_j^*(\eta)>r_j\;(j=1,2,\cdots)$
where the lower bound holds also if $j=0$. We also have the limiting values 
\begin{equation}\label{s53_limitrj*}
r_j^*(-\infty)=r_j\;(j\ge 0),\quad r_j^*(\infty)=r_{j-1}\; (j\ge 1)
\end{equation}
which are also illustrated by Figure \ref{figure:r*}. We have thus derived (\ref{s32_nueta}) and (\ref{s32_phiS*}).

We discuss the matching of the case $\rho-1=O(K^{-2/3})$ to the cases $\rho<1$ and $\rho>1$. For $\eta\to -\infty$, (\ref{s53_limitrj*}) shows that (\ref{s32_phiS*}) becomes (\ref{s32_phiS}) (for all $j\ge 0$) and the matching of the eigenvalues follows from (\ref{s53_nu}), after we replace $\widetilde{\nu}^{(0)}$ by $-r_j^*(\eta)$, and let $\eta\to -\infty$. The matching for $\eta\to +\infty$ for $j\ge 1$ then follows from the second equality in (\ref{s53_limitrj*}), but the case $j=0$ requires a separate analysis. Now $r_0^*(\eta)\to +\infty$ and we examine (\ref{s53_rj*}) with $j=0$ and $\eta\to \infty$ by using the approximation
\begin{equation*}\label{s53_Aiz}
\mathrm{Ai}(z)=\frac{1}{2\sqrt{\pi}}z^{-1/4}\exp\Big(-\frac{2}{3}z^{3/2}\Big)\big[1+O(z^{-3/2})\big],\quad z\to \infty
\end{equation*}
which leads to 
\begin{equation*}\label{s53_Ai'/Ai}
\frac{\mathrm{Ai}'(z)}{\mathrm{Ai}(z)}=-\sqrt{z}-\frac{1}{4z}+O(z^{-5/2}).
\end{equation*}
Thus, from (\ref{s53_rj*}), $\sqrt{r_0^*}+1/(4r_0^*)\sim\eta/2$ and hence 
\begin{equation}\label{s53_r0*}
r_0^*(\eta)=\frac{\eta^2}{4}-\frac{1}{\eta}+O(\eta^{-4}),\quad \eta\to +\infty.
\end{equation}
With (\ref{s53_r0*}) we see that $\nu_0$ in (\ref{s32_nu0}) behaves for $\rho\downarrow 1$ as 
\begin{equation*}
\frac{1}{K}+\frac{1}{\rho-1}\frac{1}{K^2}=\frac{1}{K}+\frac{1}{K^{5/3}}\frac{1}{\eta},
\end{equation*}
which agrees with the behavior (\ref{s53_nu}) (or (\ref{s32_nueta})) as $\eta\to\infty$.

We proceed to calculate the $O(K^{-5/3})$ correction term(s) in (\ref{s32_nueta}) explicitly. First we note that (\ref{s53_nu}) can be refined to 
\begin{equation}\label{s53_refine}
\nu=\frac{1}{K}+\frac{1}{K^{4/3}}\Big[\frac{\eta^2}{4}-r_j^*(\eta)\Big]+\frac{1}{K^{5/3}}\Big[\frac{\eta}{2}+\widetilde{\nu}^{(1)}\Big]+O(K^{-2}),
\end{equation}
where we simply expanded (\ref{s51_nuexpand}) using $\rho=1+\eta K^{-2/3}$. The correction terms $\widetilde{\nu}^{(1)}$ and $\Phi^{(1)}(S)$ again satisfy (\ref{s51_Phi1ODE}), but now $\widetilde{\nu}^{(0)}=-r_j^*(\eta)$.
By multiplying (\ref{s51_Phi1ODE}) by $\Phi^{(0)}(S)=\mathrm{Ai}(S+r_j^*)$ and integrating from $S=0$ to $S=\infty$ we obtain
\begin{equation*}\label{s53_nu1}
\widetilde{\nu}^{(1)}=\frac{-\frac{1}{2}\big[\Phi^{(0)}(0)\big]^2+\int_0^\infty S^2\big[\Phi^{(0)}(S)\big]^2dS}{\int_0^\infty\big[\Phi^{(0)}(S)\big]^2dS},
\end{equation*}
which differs from (\ref{s51_solvab}), since in the present case $\Phi^{(0)}(0)=\mathrm{Ai}(r_j^*)\ne 0$. Using the formulas in (\ref{s51_alb01})--(\ref{s51_alb03}), we obtain after some calculation
\begin{equation}\label{s53_nuk(1)}
\widetilde{\nu}^{(1)}=\widetilde{\nu}^{(1)}_j=\frac{4}{\eta^2-4r_j^*}\Big[-\frac{3}{10}-\frac{2}{15}\eta r_j^*+\frac{2}{15}\eta^2(r_j^*)^2-\frac{8}{15}(r_j^*)^3\Big],
\end{equation}
which yields, with (\ref{s53_refine}), the third term in the expansion(s) of the eigenvalues. 

The correction term $\Phi^{(1)}(S)$ in (\ref{s51_Phiexpand}) is obtained by solving the inhomogeneous Airy equation (\ref{s51_Phi1ODE}) with $\widetilde{\nu}^{(0)}=r_j^*$, $\widetilde{\nu}^{(1)}$ in (\ref{s53_nuk(1)}), and the boundary condition  
$$\frac{d}{dS}\Phi^{(1)}(S)\Big|_{S=0}+\frac{\eta}{2}\Phi^{(1)}(0)=0.$$
We omit the derivation and only give the following result:
\begin{equation}\label{s51_Phi(1)rho=1}
\Phi^{(1)}(S)=\Big(\frac{3}{10}S+\frac{19}{30}r_j^*\Big)\mathrm{Ai}(S+r_j^*)+\Big(\frac{1}{5}S^2-\frac{4r_j^*}{15}S+\frac{8}{15}(r_j^*)^2-\widetilde{\nu}_j^{(1)}\Big)\mathrm{Ai}'(S+r_j^*).
\end{equation}
This completes the analysis of the $S$ scale.

Since $\Phi^{(0)}(S)=\mathrm{Ai}(S+r_j^*)$ does not vanish as $S\to 0$, here we do not need to consider the scale $K-n=l=O(1)$; for $l=O(1)$ we would simply obtain $\phi_j(n)\sim\rho^{-n/2}\mathrm{Ai}(r_j^*)$. We do need new expansions on the $\xi$ and $n=O(1)$ scales, but that analysis is completely analogous to the case $\rho<1$. On the $\xi$-scale (\ref{s51_xiscale}) and (\ref{s51_WKB}) apply, the only difference being that the matching condition (\ref{s51_order2}) has $-r_j^*(\eta)$ in the right hand side, and so we obtain (\ref{s32_psixi*}) rather than (\ref{s32_psixi1}). For $n=O(1)$, (\ref{s32_phio1}) again applies.

\subsection{Tail behaviors of the unconditional density}\label{section:4.4}

We use (\ref{s2_cj}) and (\ref{s2_fc_ptsim}) to derive the tail behavior of the unconditional sojourn time density $p(t)$. We have, from (\ref{s2_cj}), 
$$c_0=\frac{\sum_{n=0}^{K-1}\rho^n\,\phi_0(n)}{\sum_{n=0}^{K-1}\rho^n(n+1)\,\phi_0^2(n)}\equiv \frac{\mathcal{N}}{\mathcal{D}},$$
so the tail of the unconditional density in (\ref{s2_fc_ptsim}) can be rewritten as 
\begin{equation}\label{S44_p(t)}
p(t)\sim \frac{1-\rho}{1-\rho^K}\,\frac{\mathcal{N}^2}{\mathcal{D}}\,e^{-\nu_0t}.
\end{equation}

We first consider the case $\rho<1$. The sum in $\mathcal{D}$ concentrates on the $S$-scale, and using (\ref{s32_phiS}) we have
\begin{eqnarray}\label{S44_rho<1_D}
\mathcal{D}&\sim& k_0^2 \sum_{n=0}^{K-1}(n+1)\big[\mathrm{Ai}(S+r_0)\big]^2\nonumber\\
&\sim& k_0^2\,K^{5/3}\int_0^\infty \big[\mathrm{Ai}(S+r_0)\big]^2\,dS=k_0^2\, K^{5/3}\big[\mathrm{Ai}'(r_0)\big]^2,
\end{eqnarray}
where we used the Euler--Maclaurin summation formula to approximate the sum by an integral and set $n+1=K-K^{2/3}S+1\sim K$. To approximate $\mathcal{N}$, we note that in view of the factor $\rho^n$ the sum is concentrated in the range $n=O(1)$, so we use (\ref{s32_phio1}), which yields
\begin{eqnarray}\label{S44_rho<1_N}
\mathcal{N}&\sim& k_2\sum_{n=0}^{K-1}\rho^{n/2}\frac{1}{2\pi i}\oint\frac{1}{z^{n+1}}\frac{1}{1-z}\exp\Big(\frac{1}{1-z}\Big)dz\nonumber\\
&\sim&\frac{k_2}{2\pi i}\oint\frac{1}{(1-z)(z-\sqrt{\rho})}\exp\Big(\frac{1}{1-z}\Big)dz=\frac{k_2}{1-\sqrt{\rho}}\,\exp\Big(\frac{1}{1-\sqrt{\rho}}\Big),
\end{eqnarray}
where the contour integral is a small loop about $z=0$ with $\sqrt{\rho}<|z|<1$. We also have, from (\ref{s32_k0andk1}) and (\ref{s32_k1andk2}),
\begin{equation}\label{S44_k2andk0}
k_2=k_0\frac{1}{\sqrt{e}}\,K^{1/6}\,\exp\Big\{-\frac{\pi}{2}\sqrt{K}-\frac{\pi}{4}r_0K^{1/6}\Big\}.
\end{equation}
Thus, using (\ref{S44_rho<1_D})--(\ref{S44_k2andk0}) in (\ref{S44_p(t)}), we obtain (\ref{S35_rho<1}).

Next we consider $\rho>1$. Now $\rho^n$ concentrates where $n=K-O(1)$, so we use the approximation to $\phi_0(n)$ on the $\xi$-scale in (\ref{s32_phi0xi}) with $\xi=1$, which leads to 
\begin{equation}\label{S44_rho>1_N}
\mathcal{N}\sim k_0^*\frac{\rho^K}{\rho-1}\exp\Big(-\frac{1}{\rho-1}\Big),
\end{equation}
and
\begin{equation}\label{S44_rho>1_D}
\mathcal{D}\sim (k_0^*)^2K\frac{\rho^K}{\rho-1}\exp\Big(-\frac{2}{\rho-1}\Big).
\end{equation}
Using (\ref{S44_rho>1_N}) and (\ref{S44_rho>1_D}) in (\ref{S44_p(t)}), we obtain the exponential density in (\ref{S35_rho>1}).

Now we consider $\rho\sim 1$. We first consider $\rho-1=O(K^{-2/3})$, with $\eta$ defined in (\ref{s32_eta}). Both $\mathcal{N}$ and $\mathcal{D}$ concentrate on the $S$-scale, with $\phi_0(n)$ given in (\ref{s32_phiS*}). Then we again use the Euler--Maclaurin summation formula and notice that $\rho^n\sim \exp\big(\eta K^{1/3}-\eta S\big)$, which leads to 
\begin{equation}\label{S44_eta_N}
\mathcal{N}\sim k_0\,K^{2/3}\exp\Big(\frac{\eta}{2}K^{1/3}\Big)\int_0^\infty e^{-\eta S/2}\mathrm{Ai}(S+r_0^*(\eta))\,dS,
\end{equation}
and
\begin{equation}\label{S44_eta_D}
\mathcal{D}\sim (k_0)^2 K^{5/3}\int_0^\infty \big[\mathrm{Ai}(S+r_0^*(\eta))\big]^2dS= (k_0)^2K^{5/3}\Big(\frac{\eta^2}{4}-r_0^*(\eta)\Big)\big[\mathrm{Ai}(r_0^*(\eta))\big]^2.
\end{equation}
To approximate $\frac{1-\rho}{1-\rho^K}$ we need to consider $\eta>0$ and $\eta<0$ separately, and we have
\begin{equation}\label{S44_eta_rho}
\frac{1-\rho}{1-\rho^K}\sim
\left\{ \begin{array}{ll}
\eta\,K^{-2/3}\exp\Big(-\eta K^{1/3}\Big), & \eta>0\\
-\eta\,K^{-2/3}, & \eta<0
\end{array} \right..
\end{equation}
Using (\ref{S44_eta_N})--(\ref{S44_eta_rho}) in (\ref{S44_p(t)}) yields (\ref{S35_rho1_eta-}) for $\eta<0$ and (\ref{S35_rho1_eta+}) for $\eta>0$.
Letting $\eta\to +\infty$ in (\ref{S35_rho1_eta+}) and noticing that $r_0^*(\eta)= \eta^2/4-1/\eta+O(\eta^{-4})$, we can also obtain (\ref{S35_rho>1}) as a limiting case of (\ref{S35_rho1_eta+}).

We also need to consider $\rho=1+O(K^{-1})$ and $\rho=1+O(K^{-1/2})<1$. For the $a$-scale with $\rho=1+a/K$, $-\infty<a<\infty$, we can still use (\ref{S44_eta_N}) and (\ref{S44_eta_D}), with now $\eta=aK^{-1/3}\to 0$, and we also have
$$\frac{1-\rho}{1-\rho^K}\sim\frac{a}{e^a-1}\frac{1}{K},$$
which leads to (\ref{S35_rho1_a}). 

Finally we consider the $b$-scale with $\rho=1+b/\sqrt{K}$ and $b<0$. Setting $n=K\xi$, we have $\rho^n\sim \exp\big(-b^2\xi/2+b\xi\sqrt{K}\big)$. We use the approximation to $\phi_0(n)$ in (\ref{s32_phixi}) and we have
\begin{eqnarray}\label{S44_b_N1}
\mathcal{N} &\sim& k_1\sum_{n=0}^{K-1}\rho^{n/2} K^{-1/12}\big[\xi(1-\xi)\big]^{-1/4}\exp\Big\{\sqrt{K}\psi(\xi)+K^{1/6}\psi_0^{(1)}(\xi)\Big\}\nonumber\\
&\sim& k_1K^{11/12}\int_0^1 \big[\xi(1-\xi)\big]^{-1/4}\exp\Big\{\sqrt{K}\,f(\xi)+K^{1/6}\psi_0^{(1)}(\xi)-\frac{b^2}{4}\xi\Big\}d\xi,
\end{eqnarray}
where we again used the Euler--Maclaurin summation formula and $f(\xi)$ is given by 
$$f(\xi)=\frac{b}{2}\xi+\sqrt{\xi(1-\xi)}-\sin^{-1}\sqrt{1-\xi}.$$
For $b<0$ the equation $f'(\xi)=0$ has a unique solution at $\xi=\xi_0=4/(b^2+4)$. Using the Laplace method, (\ref{S44_b_N1}) leads to 
\begin{eqnarray}\label{S44_b_N}
\mathcal{N}&\sim& k_1K^{8/12}\frac{4\sqrt{\pi}}{\sqrt{b^2+4}}\exp\bigg\{\frac{1}{b}\sin^{-1}\Big(\frac{|b|}{\sqrt{b^2+4}}\Big)-\frac{b^2+2}{b^2+4}\bigg\}\nonumber\\
&&\times\exp\bigg\{-\sin^{-1}\Big(\frac{|b|}{\sqrt{b^2+4}}\Big)\sqrt{K}+\frac{r_0}{2}K^{1/6}\Big[\frac{2b}{b^2+4}-\sin^{-1}\Big(\frac{|b|}{\sqrt{b^2+4}}\Big)\Big]\bigg\}.
\end{eqnarray}
The approximation for $\mathcal{D}$ in (\ref{S44_rho<1_D}) still holds, and for $K\to\infty$ and $b<0$ we have
\begin{equation}\label{S44_b_rho}
\frac{1-\rho}{1-\rho^K}\sim-\frac{b}{\sqrt{K}}.
\end{equation}
Using (\ref{S44_rho<1_D}), (\ref{S44_k2andk0}), (\ref{S44_b_N}) and (\ref{S44_b_rho}) in (\ref{S44_p(t)}) leads to (\ref{S35_rho1_b}). We note that letting $b\to 0^-$ in (\ref{S35_rho1_b}) we asymptotically match to the result in (\ref{S35_rho1_eta-}) as $\eta\to -\infty$. Letting $b\to -\infty$ we match to the result in (\ref{S35_rho<1}) as $\rho\uparrow 1$.

\section{Numerical studies}\label{section:5}
We assess the accuracy of our asymptotic results, and their ability to predict qualitatively and quantitatively the true eigenvalues/eigenvectors.

\begin{figure}[h]   
  \begin{minipage}[t]{0.5\linewidth}  
    \centering   
	\includegraphics[width=0.9\textwidth]{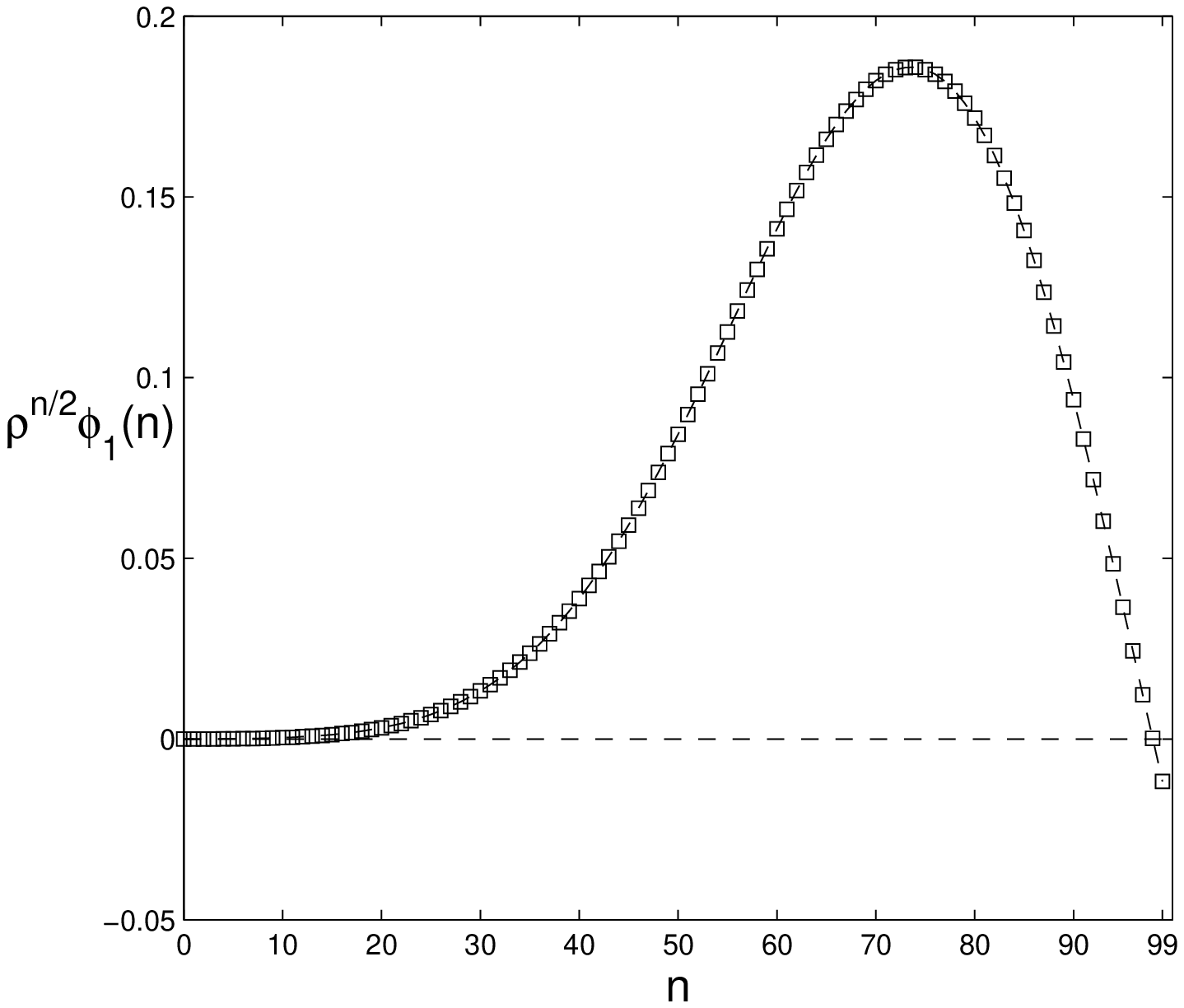}
	\caption{$\rho^{n/2}\phi_1(n)$ for $\rho=4$.} 
    \label{figure:rho>1_01}   
  \end{minipage}%
  \begin{minipage}[t]{0.5\linewidth}   
    \centering   
	\includegraphics[width=0.9\textwidth]{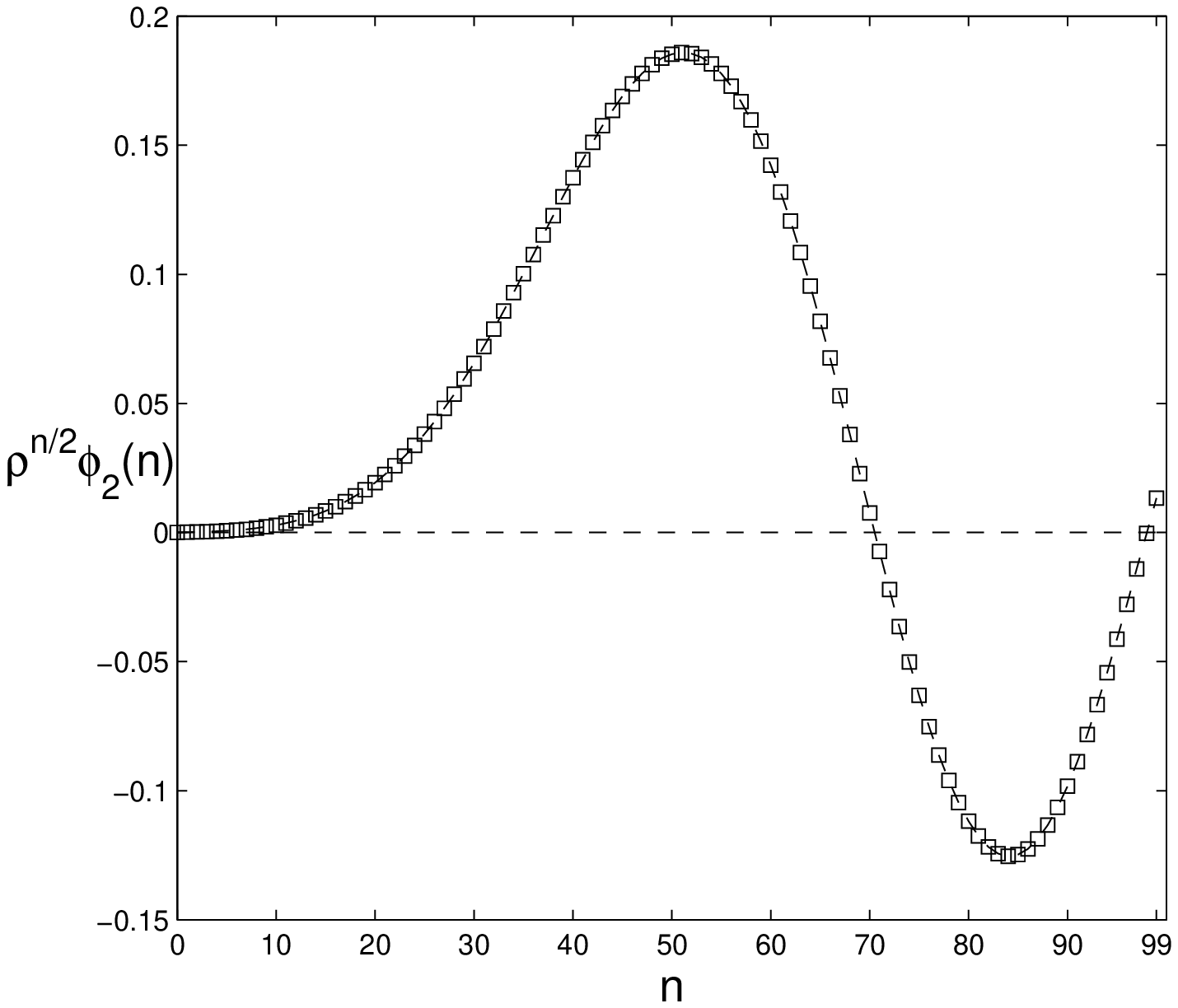}
	\caption{$\rho^{n/2}\phi_2(n)$ for $\rho=4$.} 
    \label{figure:rho>1_02} 
  \end{minipage}   
\end{figure}

\begin{figure}[h]   
  \begin{minipage}[t]{0.5\linewidth}  
    \centering   
	\includegraphics[width=0.9\textwidth]{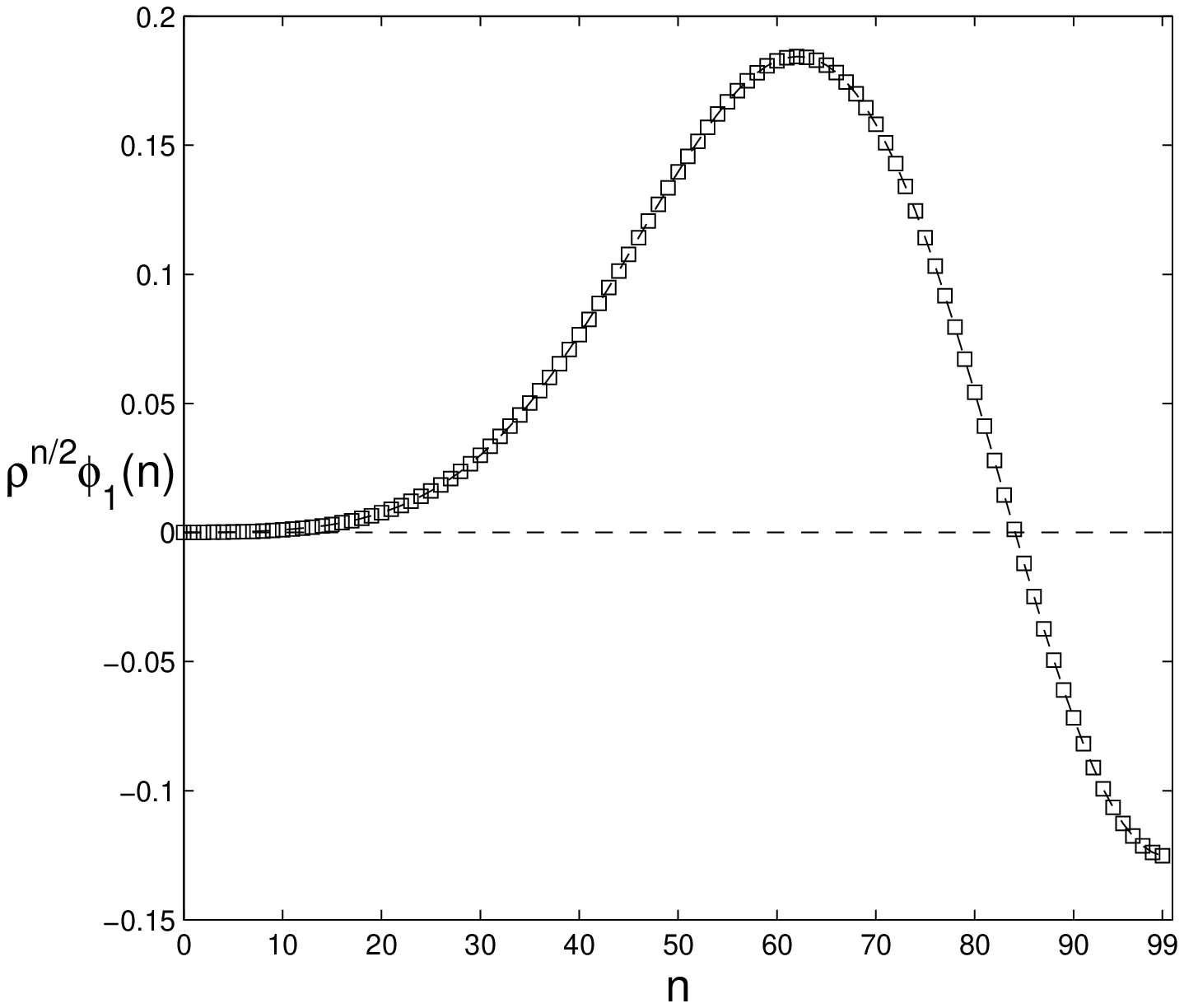}
	\caption[]{$\rho^{n/2}\phi_1(n)$ for $\rho=1$.} 
    \label{figure:rho=1_01}   
  \end{minipage}%
  \begin{minipage}[t]{0.5\linewidth}   
    \centering   
	\includegraphics[width=0.9\textwidth]{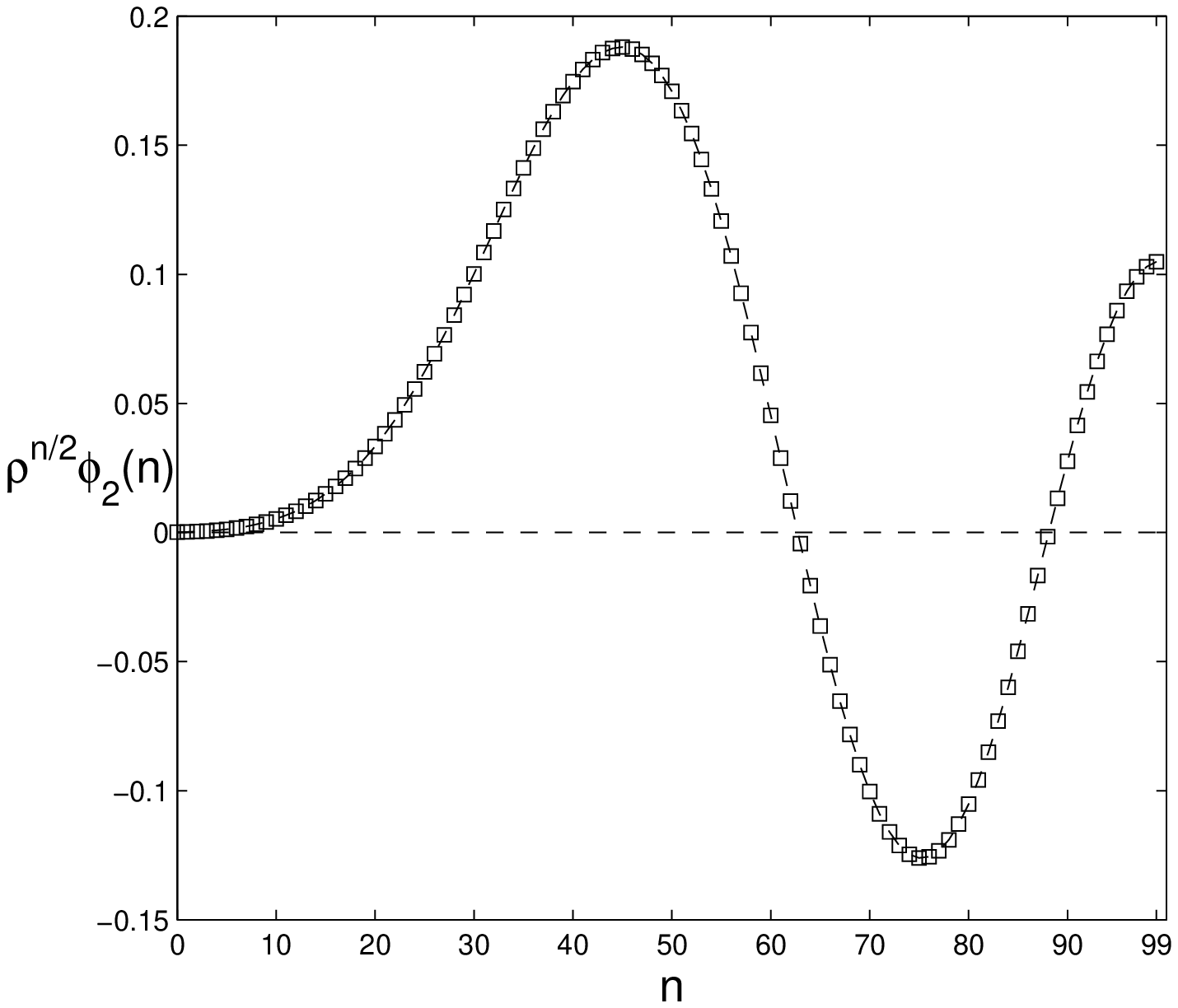}
	\caption{$\rho^{n/2}\phi_2(n)$ for $\rho=1$.} 
    \label{figure:rho=1_02}   
  \end{minipage}   
\end{figure}

\begin{figure}[h]   
  \begin{minipage}[t]{0.5\linewidth}  
    \centering   
	\includegraphics[width=0.9\textwidth]{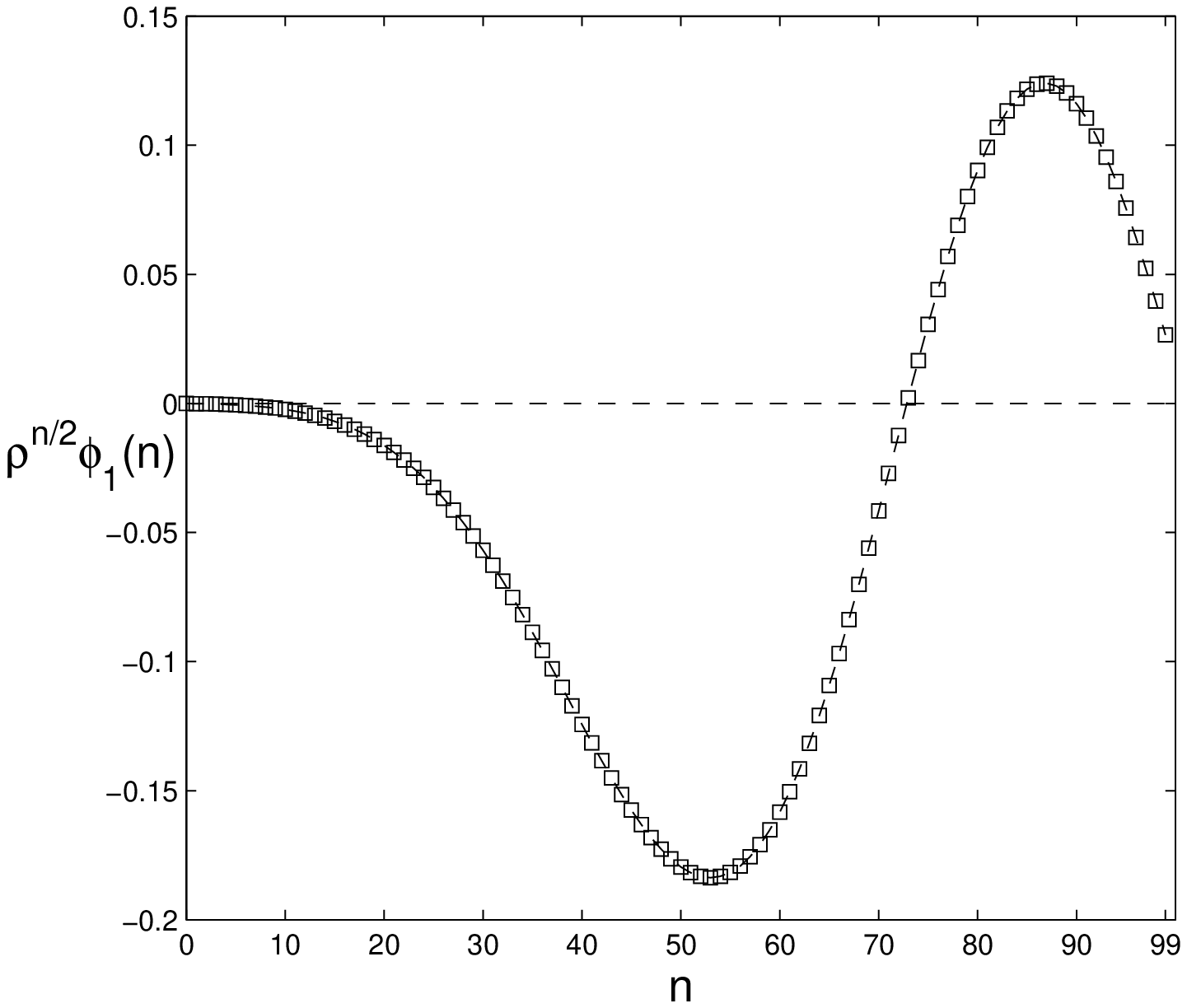}
	\caption[]{$\rho^{n/2}\phi_1(n)$ for $\rho=0.25$.} 
    \label{figure:rho<1_01}   
  \end{minipage}%
  \begin{minipage}[t]{0.5\linewidth}   
    \centering   
	\includegraphics[width=0.9\textwidth]{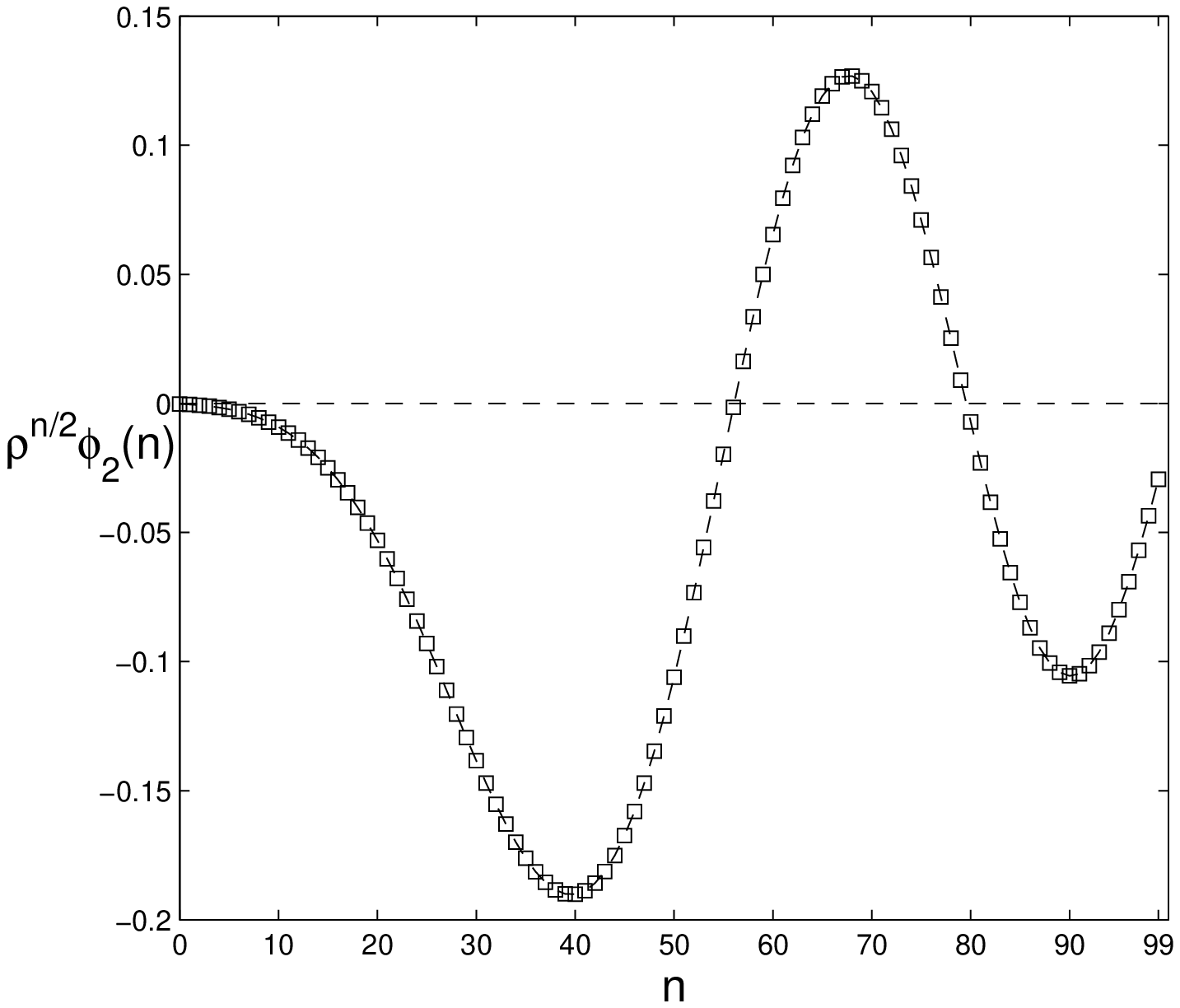}
	\caption[]{$\rho^{n/2}\phi_2(n)$ for $\rho=0.25$.} 
    \label{figure:rho<1_02}   
  \end{minipage}   
\end{figure}

In Figures \ref{figure:rho>1_01}--\ref{figure:rho<1_02} we plot the ``symmetrized" eigenvectors $\rho^{n/2}\phi_j(n)\equiv \varphi_j(n)$ for $j=1,\,2$ and the three cases of $\rho$, with $K=100$ and $n\in[0,99]$ in each case. In Figure \ref{figure:rho>1_01} we take $\rho=4>1$ and plot the first symmetrized eigenvector $\rho^{n/2}\phi_1(n)=\varphi_1(n)$. We recall that the analysis in section \ref{section:4.3} predicts that, when $\rho>1$, $\varphi_1(n)$ has a single sign change in the range $l=K-n=O(1)$. Indeed, in (\ref{s32_phil2}), when $\rho=4$, the last factor vanishes when $l=2$, is negative for $l=1$, and positive for $l\ge 3$. In Figure \ref{figure:rho>1_01}, $\varphi_1(n)$ changes sign when $n$ changes from 98 to 99. When $n=98$, $\varphi_1(n)$ is not exactly zero, but its numerical value is only about $2\%$ of the values of $\varphi_1(97)$ and $\varphi_1(99)$. Thus the asymptotic analysis gives a very good qualitative description of the true eigenvector. Note also that on the $S$-scale we have $\varphi_1(n)$ asymptotically proportional to $\mathrm{Ai}(S+r_0)$, which has a single maximum at $S\approx 1.319$, which corresponds to $n=K-K^{2/3}S\approx 71.58$, in good agreement with Figure \ref{figure:rho>1_01}, where $\varphi_1(n)$ is peaked at $n=74$. 

In Figure \ref{figure:rho>1_02} we consider $\varphi_2(n)$. For the second eigenvector the asymptotics predict one sign change on the scale $l=O(1)$ and a second on the $S$-scale. The exact sign changes occur from $n=70$ to 71, and from $n=98$ to 99. The second change is similar to that of the first eigenvector $\varphi_1(n)$, and indeed the sign change predicted by (\ref{s32_phil2}) is independent of the eigenvalue index $j$. The first sign change we predicted to occur at a root of $\mathrm{Ai}(S+r_1)$ and this function has for $S>0$ a unique root, at $S=r_0-r_1\approx 1.7498$, and this corresponds to $n=K-K^{2/3}S\approx 62.30$. To get a better approximation, we include the correction term $K^{-1/3}\Phi^{(1)}(S)$ defined in (\ref{s51_Phi(1)}). Then the two-term approximation to $\varphi_1(n)$ predicts a sign change at $n\approx 73.34$.

In Figures \ref{figure:rho=1_01}--\ref{figure:rho=1_02} we take $\rho=1$ and plot the exact $\phi_j(n)$ ($=\varphi_j(n)$). Our asymptotic results now predict that there will be $j$ sign changes on the $S$-scale, with $S>0$, and none on the $l$- and $\xi$-scales. We recall that if $\rho=1$, then $\eta=0$ in (\ref{s32_eta}), and the roots $r_j^*(0)$ of (\ref{s32_rj*}) are precisely the roots of $\mathrm{Ai}'(z)=0$; hence $r_0^*(0)\approx -1.019$, $r_1^*(0)\approx -3.248$ and $r_2^*(0)\approx -4.820$. When $j=1$, $\phi_1(n)$ in Figure \ref{figure:rho=1_01} has a single sign change as increases from $n=84$ to $n=85$, whereas the asymptotics predicts a change where $\mathrm{Ai}(S+r_1^*(0))=0$, and this occurs at $S=r_0-r_1^*(0)\approx 0.9101$, or $n\approx 80.39$. Including the correction term $K^{-1/3}\Phi^{(1)}(S)$ defined by (\ref{s51_Phi(1)rho=1}), we have $n\approx 85.86$ for the sign change. For $j=2$, $\phi_2(n)$ in Figure \ref{figure:rho=1_02} has sign changes from $n=88$ to 89, and $n=62$ to 63. The one-term asymptotic approximation in (\ref{s32_phiS*}) has zeros at $S=r_1-r_2^*(0)\approx 0.7321$ and $S=r_0-r_2^*(0)\approx 2.4820$, or $n\approx 84.23$ and $n\approx 46.53$. The two-term asymptotics approximate the sign changes as $n\approx 91.01$ and $n\approx 61.51$.

In Figures \ref{figure:rho<1_01}--\ref{figure:rho<1_02} we take $\rho=0.25<1$ and $K=100$. We again plot the symmetrized eigenvectors $\rho^{n/2}\phi_j(n)=\varphi_j(n)$. Recall that now the zeroth eigenvector $\phi_0(n)$ is given by (\ref{s32_phiS}), and has no sign changes. But, for $j\ge 1$, $\rho^{n/2}\phi_j(n)=\varphi_j(n)$ was predicted to have $j$ sign changes, in view of (\ref{s32_phiS}), and these occur at the $j$ zeros of $\mathrm{Ai}(S+r_j)$, i.e., at $S=r_k-r_j$ for $k=0,\,1,\cdots,\,j-1$. 

From Figures \ref{figure:rho<1_01} and \ref{figure:rho<1_02} we see that $\varphi_1(n)$ changes sign from $n=72$ to $n=73$ and $\varphi_2(n)$ has sign changes from $n=56$ to $n=57$, and from $n=79$ to $n=80$. In each case we also see that $\varphi_j(n)$ has no sign changes near the boundary $n=K-1(=99)$, which is consistent with the $l$-scale result in (\ref{s32_phil}) when $\rho<1$. The zero of $\mathrm{Ai}(S+r_1)$ is at $S=r_0-r_1\approx 1.7498$ or $n\approx 62.30$, while the zeros of $\mathrm{Ai}(S+r_2)$ correspond to $n\approx 31.44$ and $n\approx 69.14$. Two-term asymptotics again give better approximations, which approximate the sign change of $\varphi_1(n)$ at $n\approx 73.34$ and those of $\varphi_2(n)$ at $n\approx 48.39$ and $n\approx 80.79$.

We have thus shown that the leading-term asymptotics predict qualitatively the behavior of $\rho^{n/2}\phi_j(n)=\varphi_j(n)$ for small $j$, but the quantitative agreement is not very good when only the leading terms are used, for the moderately large value of $K=100$.
Including the correction terms to (\ref{s32_phiS}), (\ref{s32_phiS2}) and (\ref{s32_phiS*}) is beneficial.

Next we consider the accuracy of our expansions for the eigenvalues $\nu_j$. In Table \ref{table:1} we take $\rho=0.25$ and increase $K$ from 10 to 100, giving the two-term, three-term and four-term approximations to the smallest eigenvalue $\nu_0$, that result from (\ref{s32_nu}). Recall that $\nu_0$ (in fact every $\nu_j$) $\to(1-\sqrt{\rho})^2=0.25$ as $K\to\infty$. The agreement is generally very good, and shows the usefulness of the third and fourth terms in (\ref{s32_nu}). Once $K$ reaches 100, the four-term approximation agrees with the exact result to three significant figures. In Table \ref{table:2} we again increase $K$ from 10 to 100 but now take $\rho=4$. We recall that now $\nu_0=O(K^{-1})$ is given by (\ref{s32_nu0}), while the $\nu_j=O(1)$ for $j\ge 1$ follow from (\ref{s32_nuj}). The three-term approximation in (\ref{s32_nu0}) is nearly identical to the exact result, and this is consistent with the error term in (\ref{s32_nu0}) being $O(K^{-4})$. The four-term approximation to $\nu_1$ resulting from (\ref{s32_nuj}) is also very accurate, but now the error is larger, $O(K^{-2})$. 

We next show how rapidly the unconditional sojourn time density settles to its tail behavior, for moderately large $K$. We first compute $p(t)$ exactly (numerically) using (\ref{s1_pt}) and (\ref{s2_fs_pntspe}), and then compute the approximation in (\ref{s2_fc_ptsim}), which only uses the zeroth eigenvalue $\nu_0$. Then in Table \ref{table:3} we compare the exact and the approximate $-t^{-1}\log\big[p(t)\big]$ with $\rho=0.25$ and $K=10$ and 20. For $K=10$, the largest eigenvalue is $\nu_0\approx 0.3638$, which we list in the last row of the table, and the second largest eigenvalue is $\nu_1\approx 0.4858$. For $K=20$, we have $\nu_0\approx 0.3022$ and $\nu_1\approx 0.3464$. 

Table \ref{table:3} shows that both the exact and approximate values approach $\nu_0$ as $t$ increases, which coincides with our analysis, though it may take fairly large times before ultimately reaching the limit $\nu_0$. Table \ref{table:3} also shows that when $t>20$ for $K=10$, and $t>55$ for $K=20$, the relative errors are lower than 1\%, again in excellent agreement with our asymptotic analysis, since we predict that for $\rho<1$, the $j=0$ term dominates for times $t\gg O(K^{4/3})$ (see (\ref{s32_nu})), which corresponds to $t\gg 21$ for $K=10$ and $t\gg 54$ for $K=20$.

\newpage

\begin{center}
\includegraphics[width=0.6\textwidth]{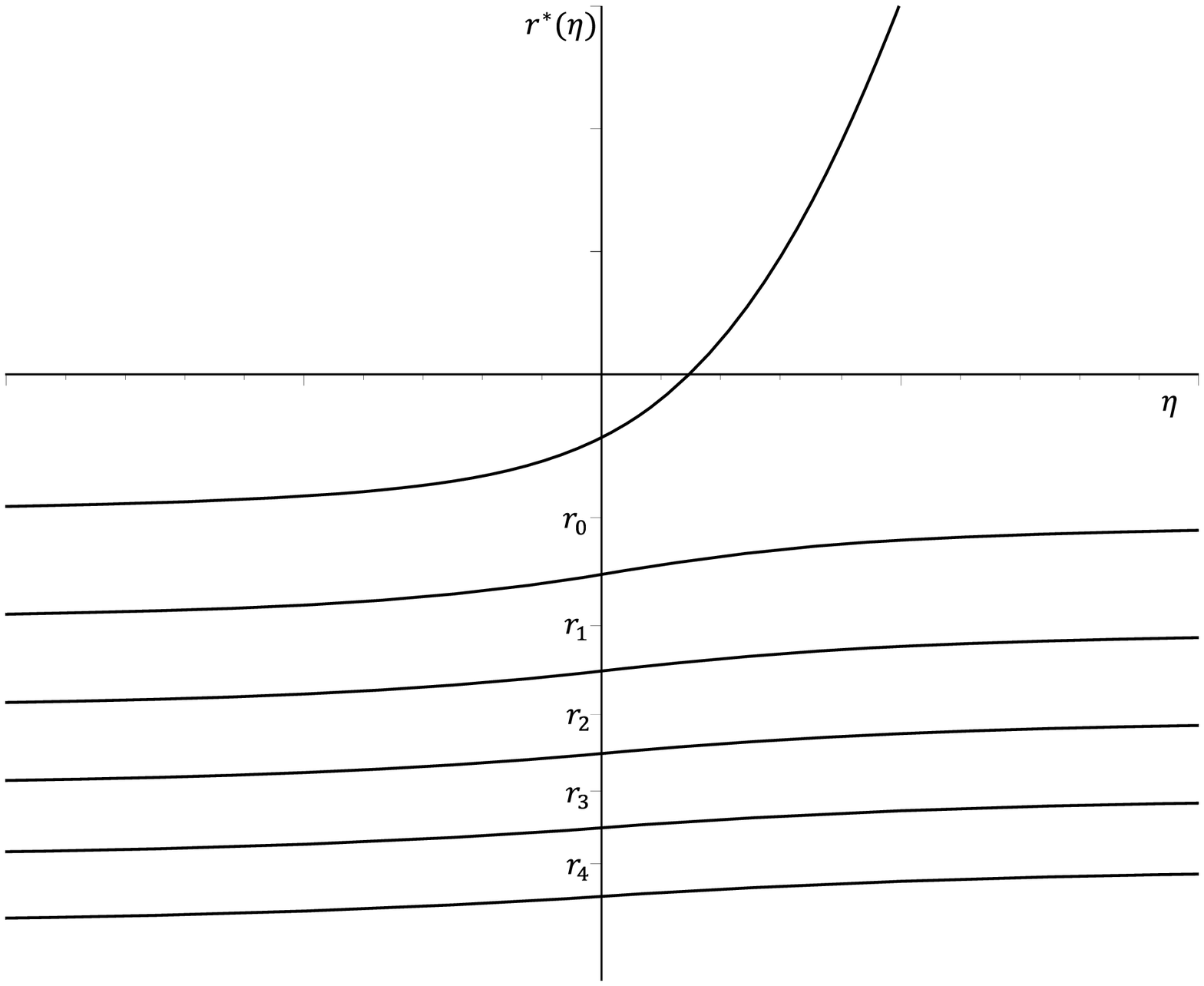}
\caption{Solution branches $r_j^*(\eta)$ of $\displaystyle \mathrm{Ai}'(r_j^*)+\frac{\eta}{2}\mathrm{Ai}(r_j^*)=0$.} \label{figure:r*}
\end{center}

\newpage

\vspace{1in}

\newpage

\begin{table}[H]
\centering
\caption{The eigenvalue $\nu_0$ with $\rho=0.25$.}\label{table:1}
  \vspace{0.5cm}
\begin{tabular} {|c||c||c|c||c|c||c|c|}       \hline
$K$		&	\tabincell{c}{$\nu_0$\\ (Exact)}		& \tabincell{c}{2-term\\ Approx.}		& \tabincell{c}{Relative\\ Error}		& \tabincell{c}{3-term\\ Approx.}		&	\tabincell{c}{Relative\\ Error}		&	\tabincell{c}{4-term\\ Approx.}		&	\tabincell{c}{Relative\\ Error}  \\\hline
10&0.3638&0.3000&1.75E-01&0.3543&2.62E-02&0.3857&6.01E-02\\\hline
30&0.2827&0.2667&5.67E-02&0.2792&1.23E-02&0.2842&5.51E-03\\\hline
50&0.2681&0.2600&3.01E-02&0.2663&6.37E-03&0.2685&1.64E-03\\\hline
70&0.2622&0.2571&1.94E-02&0.2612&3.95E-03&0.2624&7.24E-04\\\hline
100&0.2581&0.2550&1.21E-02&0.2575&2.32E-03&0.2582&3.01E-04\\\hline
\end{tabular}\\
\end{table}


\begin{table}[H]
\centering
\caption{The eigenvalues $\nu_0$ and $\nu_1$ with $\rho=4$.}\label{table:2}
  \vspace{0.5cm}
\begin{tabular} {|c||c|c|c||c|c|c|}       \hline
$K$		&	\tabincell{c}{$\nu_0$\\ (Exact)}		& \tabincell{c}{3-term \\ Approx.}		& \tabincell{c}{Relative\\ Error}		&	\tabincell{c}{$\nu_1$\\ (Exact)}		& \tabincell{c}{4-term \\ Approx.}		&	\tabincell{c}{Relative\\ Error}	  \\\hline
10&0.103413&0.103444&3.02E-04&1.634636&1.542680&5.63E-02\\\hline
30&0.033708&0.033708&8.72E-06&1.147275&1.136964&8.99E-03\\\hline
50&0.020134&0.020134&1.80E-06&1.077532&1.073979&3.30E-03\\\hline
70&0.014354&0.014354&6.42E-07&1.051445&1.049685&1.67E-03\\\hline
100&0.010033&0.010033&2.17E-07&1.033618&1.032781&8.09E-04\\\hline
\end{tabular}\\
\end{table}


\begin{table}[H]
\centering
\caption{The tail approximation of $p(t)$ with $\rho=0.25$.}\label{table:3}
\vspace{4mm}
\begin{tabular} {|c||c|c|c||c|c|c|}       \hline
&\multicolumn{3}{|c||}{$\rho=0.25,\;K=10$}  &\multicolumn{3}{c|}{$\rho=0.25,\;K=20$}  \\ \hline
$t$	& \tabincell{c}{$-\log[p(t)]/t$\\[-3pt] (Exact)} 	& \tabincell{c}{$-\log[p(t)]/t$\\[-3pt] (Approx.)}	& \tabincell{c}{Relative\\[-3pt] Error}	& \tabincell{c}{$-\log[p(t)]/t$\\[-3pt] (Exact)} 	& \tabincell{c}{$-\log[p(t)]/t$\\[-3pt] (Approx.)}		&\tabincell{c}{Relative\\[-3pt] Error} \\\hline
10	&0.5739	&0.6306	&9.88\%		&0.5734	&0.7836	&36.66\% \\\hline
20	&0.4894	&0.4972	&1.60\%		&0.4847	&0.5429	&12.02\%	\\\hline
25	&0.4671	&0.4705	&0.73\%		&0.4592	&0.4948	&7.76\%  \\\hline
40	&0.4302	&0.4305	&0.08\%		&0.4121	&0.4226	&2.53\%	\\\hline
55	&0.4123	&0.4123	&9.73E-05		&0.3860	&0.3898	&0.97\%	\\\hline
100	&0.3905	&0.3905	&2.34E-07		&0.3501	&0.3504	&0.08\%	\\\hline
1000	&0.3665	&0.3665	&$<$E-12		&0.3070	&0.3070	&$<$E-12	\\\hline
$\infty$	&0.3638	&0.3638		&--	&0.3022	&0.3022		&--\\\hline
\end{tabular}
\end{table}


\begin{center}
{\sc University of North Florida}

{\sc University of Illinois at Chicago}

\end{center}

\end{document}